
\input amstex
\documentstyle{amsppt}
\input epsf
\magnification=\magstep1
\vsize=22truecm
\newcount\escala 
\newdimen\anchura 
\newdimen\altura
\newbox\texto
\newcount\dibno                                                                 
\dibno=1 
\def\ilustracion#1,#2(#3)#4.{
\hbox to \hsize{\hfill
\altura=#1mm \anchura=#2mm
\epsfysize=\altura\epsfxsize=\anchura\epsfbox{#3.ps}\hfill}%
\nobreak{\setbox\texto=\hbox%
{\it Figure \the\dibno. #4}
\ifdim\wd\texto>\hsize{\narrower\noindent\unhcopy\texto\par}%
\else\centerline{\copy\texto}\fi}
\global\advance\dibno by 1
\ifdim\lastskip<\smallskipamount
\removelastskip\penalty55\medskip\fi
  \ignorespaces\smallbreak}

\def\ZZ{{\bold Z}}
\def\PP{{\bold P}}
\def\QQ{{\bold Q}}

\def\CC{{\bold C}}

\def\NS{\operatorname{NS}}

\def\Pic{\operatorname{Pic}}
\def\disc{\operatorname{disc}}
\def\Gal{\operatorname{Gal}}
\def\Hom{\operatorname{Hom}}
\def\rank{\operatorname{rank}}

\def\tvi{\vrule height 12truept depth 5truept width0truept}
\def\tv{\tvi\vrule}
\def\centro#1{\kern .7em\hfill#1\hfill\kern .7em}
\def\hline{\noalign{\hrule}}

\def\SS{{\Cal S}}
\def\AA{{\bold A}}

\def\EE{{\bold E}}

\def\DD{{\bold D}}
\def\BB{{\Cal B}}

\def\qed{Q.E.D.}



\newcount\parno \parno=0
\newcount\prono \prono=1
\def\sec{\S\the\parno.-\ \global\prono=1}
\def\etiqueta{\hbox{(\the\parno.\the\prono)}}
\def\finparrafo{\global\advance\parno by1
\vskip.1truecm\ignorespaces}
\def\finparrafo{\global\advance\parno by1
\vskip.1truecm\ignorespaces}
\def\cita{\ignorespaces\ \the\parno.\the\prono%
\global\advance\prono by 1}


\topmatter

\title
Miranda-Persson's problem on extremal elliptic K3 surfaces
\endtitle
\rightheadtext{Extremal elliptic $K3$ surfaces}
\author
Enrique ARTAL BARTOLO$^1$, Hiro-o TOKUNAGA$^2$ 
and De-Qi ZHANG$^3$
\endauthor
\leftheadtext{E. ARTAL, H. TOKUNAGA and ZHANG D.-Q.}
\address
Departamento de Matem\'aticas,
Universidad de Zaragoza,
Campus Plaza San Francisco s/n
E-50009 Zaragoza SPAIN
\endaddress
\email
artal\@posta.unizar.es
\endemail

\address
Department of Mathematics,
Kochi University, Kochi 780 JAPAN
\endaddress
\email
tokunaga\@math.kochi-u.ac.jp
\endemail

\address
Department of Mathematics, National University 
of Singapore, Lower Kent Rid\-ge Road, SINGAPORE 119260
\endaddress
\email
matzdq\@math.nus.edu.sg
\endemail

\thanks
$^1$ Partially supported by
CAICYT PB94-0291 and DGES PB97-0284-C02-02
\endthanks

\thanks
$^2$ Research partly
supported by the Grant-in-Aid for 
Encouragement of Young Scientists 09740031 
from the Ministry of Education, Science and Culture
\endthanks

\thanks
$^3$ Financial support
by the JSPS-NUS program
\endthanks

\endtopmatter

\document

\head\sec Introduction
\endhead

Let $f : X \to C$ be an elliptic surface over a smooth projective curve $C$ with a section $O$, i.e.,
a Jacobian elliptic fibration over $C$. Throughout this paper, we always assume that
\medbreak
($\ast$) $f$ has at least one singular fiber.
\medbreak
Let $MW(f)$ be the Mordell-Weil group of $f : X \to C$, i.e., the group of sections, $O$ being the zero. Under
the assumption ($\ast$), it is known that $MW(f)$ is a finitely generated abelian group (the Mordell-Weil
theorem). More precisely, if we let $R$ be the subgroup of the N\'eron-Severi group, $\NS(X)$, of $X$ generated 
by $O$ and all the
irreducible components in fibers of $f$, then (i) $\NS(X)$ is torsion-free, and (ii)
$MW(f) \cong \NS(X)/R$ (see [S], for instance). Note that the Shioda-Tate formula $\rank MW(f) = \rho (X) - \rank
R$ easily follows from the second statement.

We call $f: X \to C$ {\it extremal} if 

\smallbreak\item{(i)}
the Picard number $\rho (X)$ of $X$ is equal to $h^{1, 1}$ and

\smallbreak\item{(ii)}$\rank MW(f) = 0$.

\par

If $f : X \to C$ is extremal, then the Shioda-Tate formula implies $\rank R = \rho(X)$.
Hence, in other words, $f : X \to C$ is extremal 
if and only if $\rho(X) = \rank R = h^{1,1}(X)$. Also, taking the isomorphism
$MW(f) \cong NS(X)/R$ into account, it seems that we can say much about $MW(f)$ only from
the data of types of singular fibers.

\par

For extremal rational elliptic surfaces, Miranda and Persson studied them from several viewpoints \cite{MP1}; and
for such surfaces, $MW(f)$ is determined by the data of types of singular fibers. Moreover, they proved

\par
\newbox\thmp
\global\setbox\thmp=\etiqueta
\proclaim{Theorem\cita\ (\cite{MP1}, Theorem 5.4)}
For every possible configuration
of singular fibers for rational extremal elliptic surfaces, there is a unique
one with that configuration of singular fibers, except for the surfaces,
$X_{11}(j)$. These surfaces each have two singular fibers of type $I_0^*$,
with constant $J$ - map $(= j)$, and fixing $j$, there is a unique such 
surface.
\endproclaim

Suppose that $f : X \to C$ is a semi-stable elliptic $K3$ surface, i.e.,  $f$ has only $I_n$ type singular
fibers with Kodaira's notation [Ko].  In this case, $C = \PP^1$, $\NS(X) = \Pic X$, and $f$ is
extremal if and only if $f$ has exactly six singular fibers. 
For a semi-stable elliptic $K3$ surface, 
the configuration of singular fibers is said to be $[n_1,..., n_s]$ 
($n_1 \le n_2 \le \cdots \le n_s$) if it has
singular fibers $I_{n_1}$,..., $I_{n_s}$.  In \cite{MP2}, Miranda and Persson gave a complete list
for realizable $s$-tuples $[n_1,..., n_s]$; and their list shows that
there are 112 extremal cases. In \cite{MP3}, they go on to study $MW(f)$ for those
extremal elliptic $K3$ surfaces.

We say that $f : X \to \PP^1$ is of {\it type} $m$ if the corresponding 
$[n_1, n_2,..., n_6]$  
appears as the No. $m$  case in the table of \cite{MP3}.
Suppose that $f$ is of type $m$. 
What Miranda and Persson did in \cite{MP3} are that

(i) if $m \neq 2, 4, 9, 11, 13, 27, 31, 32$, 
$35, 37, 38, 44, 48, 53, 55, 69$ and $92$, 
$MW(f)$ is determined by the $6$-tuples $[n_1, n_2,..., n_6]$, and

(ii) if $MW(f) \supseteq \ZZ/2\ZZ\times \ZZ/2\ZZ$, then the corresponding elliptic
$K3$ surface is unique.

\par

The main purpose of this paper is that

\par

(i) to determine $MW(f)$ for the remaining cases, and

\par

(ii) to consider the uniqueness problem for other kinds of $MW(f)$; 
more precisely, this problem may be formulated as follows:

\par
\newbox\cuestion
\global\setbox\cuestion=\etiqueta
\proclaim{Question\cita} Let $f_1 : X_1 \to \PP^1$ and $f_2 : X_2 \to \PP^1$ be
semi-stable extremal elliptic $K3$ surfaces such that

\smallbreak\item{\rm(i)} both $X_1$ and $X_2$ have 
the same configuration of singular fibers, and

\smallbreak\item{\rm(ii)} their Mordell-Weil groups are isomorphic.

Then is it true that there exists an isomorphism $\varphi : X_1 \to X_2$
such that

\smallbreak\item{\rm(a)} $\varphi$ preserves the fibrations, and

\smallbreak\item{\rm(b)} the zero section of $f_1$ maps to that of $f_2$ with $\varphi$?
\endproclaim

\bigskip

Now let us state our result on the first problem.

\bigskip
\newbox\thmi
\global\setbox\thmi=\etiqueta
\proclaim{Theorem\cita} Let  $f : X \rightarrow \PP^1$  be
of type $m$, $m$ being one of the exceptional cases as above. 
Then we have the following table:

$$\vbox{\offinterlineskip
\halign{\tv\ #&&\tv\centro{#}&\tv#\cr
\hline
$m$ & the $6$-tuple &\hskip 5truemm $MW(f)$\hskip 5truemm & $m$ &\ the $6$-tuple
& $MW(f)$ &\cr
\hline
$2$ & $[1, 1, 1, 1, 2, 18]$ & $(0)$, $\ZZ/3\ZZ$ &
$4$ & $[1, 1, 1, 1, 4, 16]$ & $\ZZ/4\ZZ$ &\cr \hline
$9$ & $[1, 1, 1, 1, 10, 10]$ & $(0)$, $\ZZ/5\ZZ$ &
$11$ & $[1, 1, 1, 2, 3, 16]$ & $(0)$, $\ZZ/2\ZZ$ &\cr \hline
$13$ & $[1, 1, 1, 2, 5, 14]$ & $(0)$, $\ZZ/2\ZZ$ &
$27$ & $[1, 1, 1, 5, 6, 10]$ & $(0)$, $\ZZ/2\ZZ$ &\cr \hline
$31$ & $[1, 1, 2, 2, 2, 16]$ &  $\ZZ/4\ZZ$ &
$32$ & $[1, 1, 2, 2, 3, 15]$ & $(0)$, $\ZZ/3\ZZ$ &\cr \hline
$35$ & $[1, 1, 2, 2, 6, 12]$ & $\ZZ/2\ZZ$, $\ZZ/6\ZZ$ &
$37$ & $[1, 1, 2, 2, 9, 9]$ & $(0)$, $\ZZ/3\ZZ$ &\cr \hline
$38$ & $[1, 1, 2, 3, 3, 14]$ & $(0)$, $\ZZ/2\ZZ$ &
$44$ & $[1, 1, 2, 4, 4, 12]$ & $\ZZ/4\ZZ$ &\cr \hline
$48$ & $[1, 1, 2, 4, 8, 8]$ & $\ZZ/8\ZZ$ &
$53$ & $[1, 1, 3, 3, 4, 12]$ & $\ZZ/3\ZZ$, $\ZZ/6\ZZ$ &\cr \hline
$55$ & $[1, 1, 3, 3, 8, 8]$ & $(0)$, $\ZZ/2\ZZ$ &
$69$ & $[1, 2, 2, 3, 4, 12]$ &  $\ZZ/2\ZZ$, $\ZZ/4\ZZ$ &\cr \hline
$92$ & $[1, 3, 4, 4, 4, 8]$ & $\ZZ/4\ZZ$  &   &   &   &\cr \hline
}}
$$
Moreover,  all the above possibilities for  $MW(f)$  in each of
these  $17$  types, are realizable.
\endproclaim

\bigskip

Once we have settled the problem on $MW(f)$, 
we next consider Question \copy\cuestion. Our result is the following:

\medskip
\newbox\thmii
\global\setbox\thmii=\etiqueta
\proclaim{Theorem\cita} Let $f : X \to \PP^1$ 
be an extremal semi-stable elliptic $K3$ surface.
If $\sharp \left ( MW(f_i)\right ) \ge 4$, then
Question \copy\cuestion\ is true except $m = 49$ (see also Remark (0.5) (4)).
\endproclaim

\bigskip

\newbox\gnrt
\global\setbox\gnrt=\etiqueta
\remark{Remark\cita} Let $\phi$ be the homomorphism from $MW(f)$ to
 $\ZZ/n_1\ZZ\times\cdots\times\ZZ/n_6\ZZ$ given in \S 2 in \cite{MP3}, i.e.,
 $\phi(s) = (a_1,...,a_6)$, where $a_i$ is the 
 component number of the irreducible component that
 $s$ hits at the corresponding singular fiber. Since $\phi$ is injective,
 we can identify $MW(f)$ with its image by $\phi$.  Then:

\smallbreak\item{(1)}Let  $g_m : Y_m \rightarrow {\bold P}^1$
be any Jacobian elliptic fibration of type  $m$
with  $MW(g_m) = (0)$  and fitting one of the nine cases 
in Theorem \copy\thmi.  Let  $\{I_{n_1}, I_{n_2}, \dots, I_{k}, I_{k+1},$ 
$\dots, I_6\}$
be the set of types of singular fibers of  $g_m$  so that
$1 = n_1 = n_2 = \cdots = n_{k-1} < n_k \le n_{k+1} \le \cdots \le n_6$.
Then the Picard lattice  $\Pic Y_m$  is identical to
$U \oplus A_{n_k-1} \oplus \cdots \oplus A_{n_6-1}$
with the  ${\bold Q}/2{\bold Z}$-valued discriminant (quadratic) form
$q_{\Pic Y_m}$  equal to 
(cf. \cite{Mo}):
$$
(-(n_k-1)/n_k) \oplus \cdots \oplus (-(n_6-1)/n_6).
$$ 
Here  $U = \pmatrix 0&1 \\ 1&0 \endpmatrix$, and
the dual $(\Pic Y_m)^{\vee} = \Hom_{\bold Z}(\Pic Y_m, {\bold Z})$
naturally contains  $\Pic Y_m$  as a sublattice with
${\bold Z}/n_k{\bold Z} \oplus \cdots \oplus {\bold Z}/n_6{\bold Z}$
as the factor group (see \S 1 for definitions).

\smallbreak
An easy case-by-case check, using Nikulin's result that
$q_{(T_{Y_m})} = -q_{(\Pic Y_m)}$, shows that the intersection matrix
of the transcendental lattice  $T_{Y_m}$  is,
modulo the action of  $SL_2({\bold Z})$, uniquely determined 
by the data  $[n_1, \dots, n_6]$ (see \cite{Ni, Prop. 1.6.1}
or \cite{Mo, Lemma 2.4}).  So the intersection matrix
of  $T_{Y_m}$  is equal to the corresponding one in the proof of Lemma (3.3).
Thus, for each of these 9 of type  $m$,
there is exactly one $K3$ surface (modulo isomorphisms
of abstract surfaces without the fibred structure
being taken into consideration)
which has a Jacobian elliptic fibration of type  $m$
with trivial Mordell-Weil group.

\smallbreak
Also, for both  $(m, G_m) = (35, {\bold Z}/2{\bold Z}), 
(53, {\bold Z}/3{\bold Z})$, there is a unique $K3$ surface  $X_m$, 
which has a Jacobian elliptic fibration  $f_m$  of type  $m$  and  
$MW(f_m) = G_m$, because
we can prove that the transcendental lattice  $T_{X_m}$
is unique in each pair case and identical to
the corresponding one in the proof of Lemma (3.3).

\smallbreak
The authors suspect that if  $(f_m)_i : (X_m)_i \rightarrow {\bold P}^1$
are two Jacobian elliptic surfaces of the same type  $m$
and with  $MW((f_m)_1) \cong MW((f_m)_2)$ then
$(X_m)_1 \cong (X_m)_2$, though there may not be any fibred surface
isomorphism between  $((X_m)_i, (f_m)_i)$ ($i = 1,2$);
see the fourth remark below and our Proposition (4.7).  
The importance of Lemma (3.3) is that its proof can be used,
we guess, to lattice-theoretically show the existence of 
all cases of  $m$  and possibly to give an affirmative answer 
to this question.

\smallbreak\item{(2)} When  $m = 49$, we have  
$MW(f) = \ZZ/5\ZZ$
with  $s_1 = (0,0,0,2,2,2)$  or  $s_2 = (0,0,0,1,1,4)$
as its generator (cf. the Table in \cite{MP3}).  
However, we have $2s_2 = (0,0,0,2,2,10-2)$.
So we may assume that  $MW(f)$  always has  $s_1$  as its generator
after suitable relabelling of fibre components if necessary.

\smallbreak\item{(3)} When  $m = 110$, we have  
$MW(f) = \ZZ/3\ZZ \times {\ZZ}/3\ZZ$
with  
$$
G_1 = \{s_1 = (0,0,1,1,2,2), s_2 = (1,1,2,2,0,2)\}
$$
or
$$
G_2 = \{s_1 = (0,0,1,1,2,2), s_3 = (1,1,1,1,0,4)\}
$$
as its set of generators (cf. the Table in \cite{MP3}).  
Note that  $G_2$  can be replaced by the new 
generating set  $G_2' := \{s_1, 2s_3 = (3-1,3-1,2,2,0,2)\}$.
So we may assume that  $MW(f)$  always has  $G_1$  as its 
set of generators after suitable relabelling of fibre components 
if necessary.

\smallbreak\item{(4)} When  $m = 46$, we have  
$MW(f) = \ZZ/2\ZZ$
with  $s_1 = (0,0,0,0,3,5)$  or  $s_2 = (0,0,1,2,0,5)$
as its generator (cf. the Table in \cite{MP3}).  
As in the proof of Lemma (3.8), one can show that
there are pairs  $(X_i, f_i)$ ($i = 1,2$)
of the same type  $m = 46$  with
$MW(f_i) = \{O, s_i\}$.  Moreover,
the minimal resolution  $Y_i$  of  $X_i/\langle s_i \rangle$
for  $i = 1$ (resp. $i = 2$) has an elliptic
fibration  $g_i : Y_i \to {\PP}^1$,
{\it induced from}  $f_i$, of type  $m = 101$ (resp. $m = 66$).
Hence there is no isomorphism between
the pairs  $(X_i, f_i)$.

\smallbreak\item{(5)} For  $m = 69$, we have either  $MW(f) = \ZZ/2\ZZ$
with  $s = (0,1,1,0,0,6)$  as its generator, or
$MW(f) = \ZZ/4\ZZ$  with
$s = (0,1,1,0,1,3)$  as its generator (cf. Lemma (3.7).)
\endremark

\bigskip

The contents of this article is as follows:
In \S 1, we give explanations of  our technique as 
well as brief summaries of facts both of which
we need later to prove our main theorems.
In \S 2, we give a method to construct 
(or show the non-existence) of elliptic fibrations
and give several 
examples of extremal elliptic $K3$ surfaces with trivial Mordell-Weil
groups. \S 3 and \S 4 are devoted to proving Theorems \copy\thmi\
and \copy\thmii, respectively.

\remark{Acknowledgment} Part of this work was done during the second author's visit to
National University of Singapore (NUS) under the exchange program between NUS and  
Japan Society of Promotion of Science (JSPS). Deep appreciation goes to both
NUS and JSPS. The authors would like to thank Prof. S. Kondo for suggesting
Lemma (3.1).
\endremark

\definition{Conventions} In this article, 
the ground field is always the complex number field~$\CC$.

To describe the type of simple singularities 
of plane curves, we use bold capital letters, $\AA$, $\DD$ and
$\EE$.

We use capital italic letters $A$, $D$ and $E$ to describe the type of lattices, but we always multiply the
value of intersection form by $-1$ for such lattices.
\enddefinition

\finparrafo

\head\sec Preliminaries
\endhead

\subhead  1.-Cremona transformations and its applications
\endsubhead
\bigbreak

We fix notation about Cremona transformations related with
two-dimensional families of conics. Let
$V$ the vector space of homogeneous polynomials
of degree $2$ in three variables. Let $P,Q,R\in\PP^2$
three singular different points in general position
and let $V_{P,Q,R}$ be the subspace of elements of $V$
which vanish at $P$, $Q$ and $R$; it is a $3$-dimensional
vector space.
It is classical to define a rational map 
$CR_{P,Q,R}\:\PP^2\dasharrow\check\PP(V_{P,Q,R})$
where if $P_0\in\PP^2$, its image is the
hyperplane of elements of $V_{P,Q,R}$
which also vanish at $P_0$.
By a suitable choice of coordinates and the identification
of $\check\PP(V_{P,Q,R})$ with $\PP^2$ this map
may be written as:
$$
\matrix
\PP^2&\dasharrow&\PP^2\\
[x:y:z]&\mapsto&[yz:xz:xy].
\endmatrix
$$
The map $CR_{P,Q,R}$ is not defined at $P,Q,R$,
which are called the centers of the Cremona
transformation. Outside the lines joining $P,Q,R$,
this map is an isomorphism.
\medbreak

Let us consider now $P,Q\in\PP^2$ and a line $L$ through
$P$ such that $Q\notin L$. In the same way we define
$V_{P,L,Q}$ as the space of equation of conics passing
through $P$ and $Q$ and tangent to $L$ at $P$.
We define in the same way $CR_{P,L,Q}$.
We can choose coordinates such that we have:
$$
\matrix
\PP^2&\dasharrow&\PP^2\\
[x:y:z]&\mapsto&[y^2:xy:xz].
\endmatrix
$$
This map is not defined at $P$ and $Q$ and it is
an isomorphism outside $L$ and the line joining $P$ and $Q$.
We say that the centers are $P$ and the two first
infinitely near points of $P$ and $L$; we may replace
in the notation $L$ by any curve through $P$ whose
only tangent at $P$ is $L$. 
\medbreak

There is a third type of Cremona transformation
associated to a conic. Let $C$ be a smooth conic
passing through a point $P$; we denote
$V_{P,C}$ as the space of equations of conics $C'$
such that $(C\cdot C')_P=3$. We denote $CR_{P,C}$
the associated Cremona transformation. It is not
defined at $P$ and is an isomorphism outside
the tangent line to $C$ at $P$. We say that
the centers at $P$ are the three first infinitely
near points of $C$ at $P$. We can choose equations
to write it down as:
$$
\matrix
\PP^2&\dasharrow&\PP^2\\
[x:y:z]&\mapsto&[x^2:xy:y^2-xz].
\endmatrix
$$

\subhead  2.-Some lattice theory
\endsubhead
\bigbreak

We here briefly review Nikulin's lattice theory. Details are found in \cite{Ni}.
Let $L$ be a lattice, i.e., 

\smallbreak\item{(i)} $L$ is a free finite $\ZZ$ module and

\smallbreak\item{(ii)} $L$ is equipped with a 
non-degenerate bilinear symmetric pairing $\langle \,\, , \,\, \rangle$.

For a given lattice $L$, $\disc L$ is the determinant of the 
intersection matrix. Note that it
is independent of the choice of a basis.
We call $L$ unimodular if $\disc L = \pm1$.  
Let $J$ be a sublattice of $L$. We denote its orthogonal complement with respect to 
$\langle\,\, , \,\, \rangle$ by $J^{\perp}$.

For a lattice $L$, we denote its dual lattice by $L^{\vee}$. Note that, by using the pairing, $L$
is embedded in $L^{\vee}$ as a sublattice with same rank. Hence the quotient group $L^{\vee}/L$ 
is a finite abelian group, which we denote by $G_L$.

$L$ is called even if $\langle x, x \rangle$ is even for all $x \in L$. For an even lattice $L$, we define
a quadratic form $q_L$ with values in $\QQ/2\ZZ$ as follows:

$$
q_L(x \bmod L) = \langle x, x \rangle \bmod 2\ZZ.
$$

Then we have the following lemma:

\proclaim{Lemma\cita} Let $L$ be an unimodular lattice. 
Let $J_1$ and $J_2$ be  sublattices  of $L$
such that  $J_1^{\perp} = J_2$ and $J_2^{\perp} = J_1$. Then

$(i) \quad G_{J_1} \cong G_{J_2}$ and $(ii) \quad q_{J_1} = - q_{J_2}$.
\endproclaim

For a proof, see \cite{Ni}.

A sublattice $M$ of $L$ is called primitive if $L/M$ is torsion-free.

\example{Example\cita} For a $K3$ surface $X$, $H^2(X, \ZZ)$ is an even unimodular lattice with respect to
the intersection pairing. The Picard group, $\Pic X$, is a primitive sublattice of $H^2(X, \ZZ)$, and
$T_X := \left (\Pic X \right )^{\perp}$ is called the transcendental lattice of $X$.
\endexample

We shall end this subsection with the following lemma.

\newbox\lemmasi
\global\setbox\lemmasi=\etiqueta
\proclaim{Lemma\cita} For  $j = 1,2$, let  $\Delta_j = 
\Delta(1)_j \oplus \cdots \oplus \Delta(r_j)_j$
be a lattice where each  $\Delta(i)_j$  is of Dynkin type  $A_a, D_d$
or  $E_e$.

\smallbreak\item{\rm(1)} Suppose that  $\Phi : \Delta_1 \rightarrow \Delta_2$
is a lattice-isometry.  Then  $r_1 = r_2$  and
$\Phi(\Delta(i)_1) = \Delta(i)_2$  after relabelling.

\smallbreak\item{\rm(2)} Let  ${\bold A} = A_{m_1} \oplus \cdots \oplus A_{m_k}$
be a direct sum of lattices of Dynkin type  $A_{m_i}$.
Suppose that  ${\bold A}$  is an index-$n$ ($n > 1$) sublattice  
of  $\Delta := \Delta_2$  and that  $(m_1, \dots, m_k) = 
(1, 1, 5, 11), (2, 2, 3, 11)$.  Then one of the following
three cases occurs (the first two are quite unlikely but the authors
do not have a proof yet):

\smallbreak\itemitem{\rm(2-1)} ${\bold A} = 
A_1 \oplus (A_1 \oplus A_5 \oplus A_{11}),
\Delta = A_1 \oplus D_{17}$, and
$(A_1 \oplus A_5 \oplus A_{11}) \subseteq D_{17}$  is
an index-6 extension.

\smallbreak\itemitem{\rm(2-2)} ${\bold A} = 
A_2 \oplus (A_2 \oplus A_3 \oplus A_{11}),
\Delta = A_2 \oplus D_{16}$, and
$(A_2 \oplus A_3 \oplus A_{11}) \subseteq D_{16}$  is
an index-6 extension.

\smallbreak\itemitem{\rm(2-3)} ${\bold A} =
A_1 \oplus A_{11} \oplus (A_1 \oplus A_5),
\Delta = A_1 \oplus A_{11} \oplus E_6$, and
$(A_1 \oplus A_5) \subseteq E_6$  is
an index-2 extension.

\endproclaim

\demo{Proof} We observe that
$$
|\det(A_n)| = n+1,\ |\det(D_n)| = 4,\ |\det(E_5)| = 3, |\det(E_7)| = 2,
|\det(E_8)| = 1.
$$  
We also note that for an index $n$ lattice extension  
$L \subseteq M$  one has  $|\det(L)| = n^2 |\det(M)|$.

\par
(1) is true when  $r_1 = r_2 = 1$.  In general,
for a generating root  $e$  in  $\Delta(1)_1$
with  $e^2 = -2$, one has  $(\Phi(e))^2 = -2$  and hence
$\Phi(e) \in \Delta(1)_2$  say, because  $\Delta_2$
is even and negative definite.  Now the connectedness of 
$\Delta(1)_1$  implies that  $\Phi(\Delta(1)_1) \subseteq \Delta(1)_2$.
Thus to prove (1), we may assume that  $r_2 = 1, \Delta_2 = \Delta(1)_2$.
The same argument applied to  $\Phi^{-1}$  shows that
$r_1 = 1$.

\par
(2) The argument in (1) applied to the inclusion
${\bold A} \hookrightarrow \Delta_2$, implies that
each  $\Delta(i)_1$  contains a finite-index sublattice
which is a sum of a few summands of  ${\bold A}$.
Now it follows from the observations at the beginning of 
the proof of this lemma, that either (2) is true or
one of the following two cases occurs:

\par
Case (2-4) ${\bold A} = A_{11} \oplus (A_2 \oplus A_2 \oplus A_3),
\Delta = A_{11} \oplus D_7$, and
$(A_2 \oplus A_2 \oplus A_3) \subseteq D_7$  is
an index-3 extension.

\par
Case (2-5) ${\bold A} = A_2 \oplus A_3 \oplus (A_2 \oplus A_{11}),
\Delta = A_2 \oplus A_3 \oplus D_{13}$, and
$(A_2 \oplus A_{11}) \subseteq D_{13}$  is
an index-3 extension.

\par
In the following, if  $e_i$'s form a canonical 
${\bold Z}$-basis of  $A_n$  we let
$h_n = (1/(n+1))\sum_{i=1}^n i e_i$ (mod $A_n$)  be 
the generator of  $(A_n)^{\vee}/A_n \cong {\bold Z}/(n+1){\bold Z}$.  
Note that  $(h_n)^2 = -n/(n+1)$.

\par
Suppose the contrary that Case (2-4) occurs.  Set  
${\bold B} = A_2 \oplus A_2 \oplus A_3$.  
Then  $D_7 \subseteq {\bold B}^{\vee} := \Hom_{\bold Z}({\bold B}, {\bold Z})$.
and the latter is generated by  $h_2, h_2', h_3$
with  $(h_2)^2 = -2/3 = (h_2')^2, (h_3)^2 = -3/4$.
Since  $D_7$  is generated by roots and contains
an index-3 sublattice  ${\bold B}$,
there is a root  $t \in D_7 - {\bold B}$,
and we can write  $t = ah_2 + bh_2' + A$  where
$a, b \in {\bold Z}, A \in {\bold B}$.  Then  
$-2 = t^2 = (-2/3)(a^2 + b^2) + A^2 - 2s_1$  for some  
$s_1 \in {\bold Z}$.  Since  ${\bold B}$  is even and
negative definite, $A^2 = -2s_2$  for some  $s_2 \in {\bold Z}$.
Denote by  $s = s_1 + s_2$.
Then  $3 = a^2 + b^2 + 3s$, $3 | (a^2+b^2)$.  
Hence  $a = 3a_1, b = 3b_1$  for some  $a_1, b_1 \in {\bold Z}$.  
This leads to that
$t = a_1(3h_2) + b_1(3h_2') + A \in {\bold B}$, a contradiction. 

\par 
Suppose the contrary that Case (2-5) occurs.  Set
${\bold B} = A_2 \oplus A_{11}$.  
Then  $D_{13} \subseteq {\bold B}^{\vee}$
and the latter is generated by  $h_2, h_{11}$.
As in Case (2-4), there is a root  $t \in D_{13} - {\bold B}$,
and we can write  $t = ah_2 + 4bh_{11} + A$  where
$a, b \in {\bold Z}, A \in {\bold B}$.  Then
$-2 = t^2 = (-2/3)(a^2 + 22b^2) -2 s$  for some
$s \in {\bold Z}$.  
Hence  $3 = a^2 + 22b^2 + 3s$, $3 | (a^2+b^2)$
and  $a = 3a_1, b = 3b_1$
for some  $a_1, b_1 \in {\bold Z}$.  This leads to that
$t \in {\bold B}$, a contradiction.
\qed\enddemo

\subhead  3.-Review on elliptic surfaces with many torsions
\endsubhead

\bigbreak

We here give a brief summary on the results in \cite{CP} and 
\cite{C}. Let $f : X \to C$ be an elliptic surface
over a curve $C$ with a section $s_0$. Let $MW(f)$ be its Mordell-Weil group, the group of
sections, $s_0$ being the zero element. We denote its torsion part by $MW(f)_{tor}$. 
Suppose that $MW(f)_{tor} \supset \ZZ/m\ZZ\oplus\ZZ/n\ZZ,\, m|n,\,  mn \ge 3$. Then it is known that one obtains
$f : X \to C$ in a certain universal way, which we describe below. For that purpose, we need
some notations.

Set 
$$
\Gamma_m(n)=
\left\{
\pmatrix
a & b \\
c & d 
\endpmatrix
\in SL(2, \ZZ)\ \mid\ 
\pmatrix
a & b \\
c & d 
\endpmatrix 
\equiv
\pmatrix
1 & \ast \\
0 & 1 
\endpmatrix
\bmod n, \, b \equiv 0 \bmod m 
\right \}
$$
Let $X_m(n) = \Gamma_m(n)\backslash {\Cal H}^*$,
where ${\Cal H}^*$ is the upper halfplane in $\CC$, 
and let $E_m(n)$ be the elliptic modular
surface of $\Gamma_m(n)$. By definition, $E_m(n)$ is an elliptic surface over $X_m(n)$; and we denote
the morphism from $E_m(n)$ to $X_m(n)$ by $\psi_{m, n}$.

Suppose that $MW(f)_{tor} \supset \ZZ/m\ZZ\oplus\ZZ/n\ZZ,\, m|n,\,  mn \ge 3$. Then we have a commutative diagram
$$
\matrix
C&\overset g\to\rightarrow&X_1(N)\\
&j\searrow&\downarrow j_{m,n}\\
&&\PP^1
\endmatrix
$$
where $j$ and $j_{m,n}$ are the $j$-invariants of $f$ and $\psi_{m,n}$, 
respectively.
Moreover, this diagram essentially gives $f: X \to C$, i.e., $X$ is obtained as the
pull-back surface by $g$, in the sense of relatively minimal smooth model.

Thus $f$ is determined by $g$. Hence the uniqueness of $X$ is reduced to that of $g$, which we consider in 
\S 4.

\subhead 4.- Comments on pencil of plane curves and nodal cubics
\endsubhead
\bigbreak

Let $C=\{f=0\}$ and $D=\{g=0\}$ two projective plane curves
of degree $d$ without common components.
They  define a pencil of curves by
considering $\{C_{[t:s]}\}_{[t:s]\in\PP^1}$, where
$C_{[t:s]}$ is the curve of equation $sf-tg=0$.
Let us denote $\BB:=C\cap D$; it is the
set of base points of the pencils;
these base points are the intersection points of
any couple of element of the pencil. A base
point $P$ is multiple if $(C\cdot D)_P>1$
(we may replace $C$ and $D$ by any couple
of different elements of the pencil). A pencil defines
a rational map $\PP^2\dasharrow\PP^1$ which
is well-defined outside the base points. Let $Z\subset\PP^2$
be an irreducible curve of degree $e$
which is not a component of
any element in the pencil. Let $C_{[t:s]}$
a generic element of the pencil. Then the pencil
defines a map $\phi\:Z\to\PP^1$  of degree
$$
d_Z:=de-\sum_{P\in\BB}(Z\cdot C_{[t:s]})_P;
$$
if a base point $P$ is in $Z$ its image is the
unique value $\phi(P)$ such that $(Z\cdot C_{\phi(P)})_P$
is greater than the generic intersection number.
The critical points of the map are
the points $Q\in Z$ such that:

\smallbreak\item{--} If $Q$ is not a base point,
then $C_{\phi(Q)}$ is either singular at $Q$
or not transversal to $Z$ at $Q$, i.e.,
$(Z\cdot C_{\phi(Q)})_Q>1$.

\smallbreak\item{--} If $Q\in\BB$,
then $(Z\cdot C_{\phi(Q)})_Q>1+(Z\cdot C_{[t:s]})_P$,
for $[t:s]\neq\phi(Q)$.

\bigbreak
Let us consider a nodal cubic $N$ in $\PP^2$. We will
apply later the next well-known result.

\newbox\nodal
\global\setbox\nodal=\etiqueta
\proclaim{Proposition\cita} There exists a 
homogeneous coordinate system $[x:y:z]$ in $\PP^2$
such that the equation of $N$ is $xyz+x^3-y^3=0$.
The subgroup $G$ of $PGL(3,\CC)$ fixing $N$ is isomorphic
to the dihedral group of order $6$. Let $\varphi\:\CC^*\to Reg(N)$
be the mapping defining by $\varphi(t):=[t:t^2:t^3-1]$. Let us consider
on $N$ the geometrical group structure with zero element $[1:1:0]=\varphi(1)$.
Then $\varphi$ is a group isomorphism. Each element of $G$ is determined
by its action on $Reg(N)$; the induced action on $\CC^*$ is generated
by $t\mapsto t^{-1}$ and $t\mapsto\zeta t$ where $\zeta^3=1$.
\endproclaim

\finparrafo

\head\sec
Some extremal elliptic 
$K3$ surfaces with trivial Mordell-Weil group
\endhead

\subhead 1.- Elliptic fibrations and sextic curves
\endsubhead
\bigbreak
Relationship between extremal elliptic fibrations and 
maximizing sextic curves
was intensively studied in Persson's paper \cite{P}.
We explain in this section how to apply this method to
construct or discard extremal elliptic fibrations.
Let $(X,f)$ be a pair
such that $X$ is a $K3$ surface
and $f \:X\to\PP^1$ is a relatively minimal elliptic fibration
with a fixed section $O$.

\definition{Step 1} Fix $O$ as the zero element of the
Mordell-Weil group $MW(f)$. It
determines a group law on each regular fiber and it
extends to a group law in the regular part of any fiber. For a fiber
$F$ of type $I_n$, there is a short exact sequence
$$
0\to\CC^*\to Reg(F)\to \ZZ/n\ZZ\to 0
$$
where the kernel corresponds to the part of 
$Reg(F)$ in the irreducible
component which intersects $O$.
\enddefinition

\definition{Step 2} On the regular
part of any fiber $F$ we can consider the map
$P\mapsto -P$, (where $F\cap O$ is the
zero element). These maps are the
restriction of a morphism $\sigma\:X\to X$,
which is clearly an involution. By definition $f \circ\sigma= f$.
Then, there is a natural map $\tilde\rho\:X/\sigma\to\PP^1$; if
$F$ is an elliptic fiber of $\pi$, $\tilde\rho(F)$ 
is the quotient of an elliptic curve by an involution with
four fixed points (the $2$-torsion), i.e., a smooth rational curve.

Then $\tilde\rho\:X/\sigma\to\PP^1$ is a morphism from
a smooth (rational) surface onto $\PP^1$ whose generic fiber is $\PP^1$.
If $F$ is a fiber of type $I_{2n+1}$ (resp. $I_{2n}$), $\tilde\rho(F)$ 
is  a curve with normal crossings and
$n+1$ irreducible components which are smooth and
rational.
\enddefinition

\definition{Step 3} For any singular fiber $F$,
we contract
all of the irreducible components of $\tilde\rho(F)$ but the one which intersects
the $\tilde\rho(O)$. 
We obtain a holomorphic fiber bundle
$\rho\:\Sigma\to\PP^1$ with fiber isomorphic to $\PP^1$ ($\Sigma$ smooth)
and a map $\tau\:X\to\Sigma$ such that $\rho\circ\tau=\pi$.
This map is generically $2:1$.
\enddefinition

The map $\tau$ is a $2$-fold covering ramified 
on the image of the fixed points of $\sigma$, i.e.,
on the
image of the $2$-torsion. 
We can write this curve as $E\cup R$ where 
$E:=\tau(O)$, $R\cap E=\emptyset$ and
$R$ has intersection number three with the fibers
of $\rho$.
The number of irreducible components
of $R$  depends on the $2$-torsion $T_2(MW(f))$
of the Mordell-Weil group of $X$ (one irreducible
component if $T_2(MW(f))=0$, two if $T_2(MW(f))=\ZZ/2\ZZ$
and three if $T_2(MW(f))=\ZZ/2\ZZ\oplus\ZZ/2\ZZ$).
\medbreak

If the configuration of $\pi$ is $[1,\dots,n_1,\dots,n_r]$,
$1<n_1\leq\dots\leq n_r$,
then $R$ has exactly $r$ singular points of type
$\AA_{n_1-1},\dots,\AA_{n_r-1}$.

\remark{Remark\cita} Let us suppose that $n_r>7$,
and let us call $F$ the fiber of $\rho$ containing
this point $\AA_{n_r-1}$; $R$ intersects also $F$ at another
point $P$. Then we can perform three
Nagata elementary transformations on the first three
infinitely near points of $R$ at $\AA_{n_r-1}$.
We call $\Sigma'$ the result of this operation
and we do not change the notation for the strict transforms;
it induces a new fibration $\rho'\:\Sigma'\to\PP^1$
where $E$ is a section
of self-intersection $-1$.
The curve $R$ has a singular point $\AA_{n_r-7}$
and $(R\cdot E)_P=3$, and $R$ is smooth at $P$. 
We can contract $E$ and we obtain a projective plane
where the contraction of $R$ is a curve
of degree $6$
(also denoted by $R$) which has $r+1$ singular points
of type $\AA_{n_1-1},\AA_{n_2-1},\dots,\AA_{n_r-7}$
and $\EE_6$; the image of $F$ is the tangent line to $R$ at $\EE_6$
and passes through $\AA_{n_1-7}$. The pencil which induces the
elliptic fibration (the {\it preferred pencil\/}) is the
pencil of lines through $\EE_6$. This fibration is called
the standard fibration in \cite{P} and in this case
$\EE_6$ is its center.
\endremark

We can consider some kind of converse of this construction. Let
$R\subset\PP^2$ be a reduced curve (maybe reducible) of degree six
such that its singular points are simple. Let $P$ be a singular point
of $R$. Then if $X$ is the minimal resolution of the
ramified double covering of $\PP^2$ ramified
on $R$ and $\pi\:X\to\PP^1$ is the mapping induced by the pencil of
lines through $P$, then $\pi$ is a relatively minimal elliptic
fibration of the $K3$-surface $X$. We call $(X,\pi)$ the elliptic
fibration associated to $(R,P)$
and we will call the pencil of lines at $P$ the preferred
pencil; we will denote
$\sigma\:X\to\PP^2$ the double covering. Next result is easy and useful.

\newbox\cover
\global\setbox\cover=\etiqueta
\proclaim{Proposition\cita} Let $\pi\:X\to\PP^1$ the elliptic
fibration associated to $(R,P)$ as above. Let $E$ be a section
of $X$; let $C:=\sigma(E)$. Then either $C$ is an irreducible component
of $R$, either the intersection number of $C$ and $E$ at any point
in $C\cap R$ is an even number.

In both cases $C$ is a curve of degree $d$ having at $P$ a singular point
of multiplicity $d-1$. In the first case there is exactly one section over
$C$ and in the second case there are exactly two such sections.
\endproclaim

We study now the existence of 
elliptic fibrations with trivial Mordell-Weil group
in the cases of ambiguity which appear in the list
of Miranda and Persson. In fact, we have applied this method 
to all cases of ambiguity in the list. As it is very long, we present
only a few cases, where interesting phenomena occur.

\subhead 2.-Type $m=9$
\endsubhead
\bigbreak

\proclaim{Proposition\cita} There exist elliptic
$K3$ surfaces of type $9$, i.e., with configuration 
$[1,1,1,1,10,10]$, and trivial Mordell-Weil group.
\endproclaim

This proposition gives one ambiguity case as such a fibration
with Mordell-Weil group of order $5$ appears in \cite{MP3}.

We look for an irreducible curve $R$ of degree 6 having
three singular points of type $\EE_6,\AA_3,\AA_9$ and such
that the tangent line to $R$ at $\EE_6$ passes through $\AA_3$.
As in the case above the line through $\AA_3$ and $\AA_9$
intersects $R$ at two other points. 

\proclaim{Step 1}First Cremona transformation.
\endproclaim

We consider $CR_{\EE_6,\AA_3,\AA_9}$. 
We denote $R_1$ the strict transform of $R$;
$R_1$ is a quintic curve. We have a smooth point $Q$
such that the tangent line $T$ to $R_1$ at $Q$
verifies that $(R_1\cdot Q)_Q=4$. We denote
$Q'$ the other point in $R_1\cap T$.

The other singular points of $R_1$ are
$\AA_7$ (coming from $\AA_9$), $P_1$
(an ordinary double point coming from $\AA_3$) and another ordinary
double point denote $P_2$.
The preferred pencil of lines has its center at $P_1$.
The line joining $P_1$ and $P_2$
intersects $R_1$ at
$Q$. The line joining $P_1$ and $\AA_7$ passes through $Q'$. 
The ramification locus is $R_1\cup T$.

\ilustracion 40,120(cremona4).

\proclaim{Step 2}Second and third Cremona transformations.
\endproclaim

We perform $CR_{P_1,P_2,\AA_7}$.
We obtain a quartic curve $R_2$ with one
singular point $\AA_5$ (coming from $\AA_7$).
The line $T$ becomes a conic $T_2$ and 
$R_2\cap T_2=\{Q,Q',Q''\}$ where
$(R_2\cdot T_2)_{Q}=5$, $(R_2\cdot T_2)_{Q'}=2$,
$(R_2\cdot T_2)_{Q''}=1$, and $\AA_5,Q',Q''$
are aligned. The center of the preferred pencil
is $Q''$.

We perform the third Cremona transformation
$CR_{\AA_5, L, Q''}$, $L$ being the tangent line at $\AA_5$.
We obtain  two cubics
$R_3$ and $T_3$. The cubic $R_3$ has an
ordinary double point $\AA_1$ and  $T_3$
has also a double point denoted $S$ (which
is the center of the preferred pencil).
The curves $R_3$ and $T_3$
have two intersection points $Q$ and $Q'$,
with intersection numbers $5$ and $4$,
and the points $Q'$, $S$ and $\AA_1$ are aligned.

\newbox\curvb
\global\setbox\curvb=\etiqueta
\proclaim{Question\cita} Do there exist an irreducible
nodal cubic $R_3$ (with node $\AA_1$), an irreducible
cubic $T_3$ with a double point $S$ in $\PP^2$ 
such that $R_3\cap T_3=\{Q,Q'\}$, $Q,Q'\neq S,\AA_1$,
with $(R_3\cdot T_3)_Q=5$, $(R_3\cdot T_3)_{Q'}=4$
and $Q',S,\AA_1$ aligned?
\endproclaim

\proclaim{Proposition\cita}
The answer to Question \copy\curvb\ is yes.
\endproclaim

\demo{Proof}
We proceed by applying
Proposition \copy\nodal\ to $R_3$.
We suppose that $Q=p(s^{-4})$ and
$Q'=p(s^5)$. In this situation the
equation of the line joining $Q'$ and $\AA_1$
is $y=s^5x$. Let $f(x,y,z)=0$ an equation
for $T_3$ such that
the coefficient of $z^3$ in $f$ is $1$.
Then $f(t,t^2,t^3-1)=(t-s^5)^4(t-s^{-4})^5$.
We impose that $T_3$ intersects
the line $y=s^2x$ at one point outside $Q'$
(with multiplicity $2$). We force this point
to be singular and we get the conditions on
$s$ (again with MapleV).
We obtain that
$$
(s^6-1)(s^6+3s^3+1)(s^{12}+4s^9+s^6+4s^3+1)=0.
$$
We consider the action of
the dihedral group; in the first term it is enough
to retain the cases $s=\pm1$; the positive case
is too degenerated so it remains only $s=-1$.
The equation of $T_3$ in this case is:
$$
13\,{y}^{3}+9\,{y}^{2}x-5\,{y}^{2}z-9\,y{x}^{2}-6\,yxz-y{z}^{2}-13\,{x
}^{3}-5\,{x}^{2}z+x{z}^{2}+{z}^{3}=0.
$$
For the second term, one can see that we
force $S=\AA_1$ which is also too degenerated.
The last factor gives two different cases
(the twelve roots give two orbits
by the action of the dihedral group).
The equation is:
$$
\multline
\left (-{\frac {1265\,{s}^{9}}{2}}-60\,{s}^{3}-{\frac {4671}{2}}-2170
\,{s}^{6}\right ){x}^{3}+\left (1205\,{s}^{8}+320\,{s}^{11}+1285\,{s}^
{2}\right )z{x}^{2}\\
+\left (10080\,s+135\,{s}^{4}+9480\,{s}^{7}+2466\,{
s}^{10}\right )y{x}^{2}+\left (60\,s+60\,{s}^{7}+16\,{s}^{10}+5\,{s}^{
4}\right ){z}^{2}x\\
+\left (15255\,{s}^{2}+216\,{s}^{5}+14325\,{s}^{8}+
3780\,{s}^{11}\right ){y}^{2}x+\left ({\frac {495\,{s}^{9}}{2}}+{
\frac {2103}{2}}+990\,{s}^{6}\right )yzx\\
+\left (-{\frac {1735\,{s}^{9}
}{2}}-60\,{s}^{3}-{\frac {6609}{2}}-3110\,{s}^{6}\right ){y}^{3}
-\left (640\,s+620\,{s}^{7}+160\,{s}^{10}+5\,{s}^{4}\right )z{y}^{2}\\
+\left (-75\,{s}^{2}-75\,{s}^{8}-20\,{s}^{11}-4\,{s}^{5}\right ){z}^{2}
y+{z}^{3}=0
\endmultline
$$
\qed
\enddemo

We deduce that there are essentially three different
answers to Question \copy\curvb.
The main feature of the first answer is that the
tangent line $L$ to $R_3$ at $Q'$ passes through $Q$.
The elliptic surface is obtained from
the double covering of $\PP^2$ ramified along
$R_3+T_3$, and the elliptic fibration comes
from the pencil of lines with center at $S$.
One of the singular fibers is produced by the
line joining $S$, $\AA_1$ and $Q'$.

\ilustracion 40,120(fibra1).

The other singular fiber is produced
by the line joining $S$ and $Q$.

\ilustracion 40,120(fibra2).

\proclaim{Proposition\cita} The solution for $s=-1$ produces
the elliptic fibration such that $MW$ is cyclic of order $5$.
The solutions $s^{12}+4s^9+s^6+4s^3+1=0$ produce
elliptic fibrations with trivial Mordell-Weil group;
this case was not previously known.
\endproclaim

\demo{Proof} We note that the exceptional curve of
the blowing-up of $S$ never produces a section.
In both cases the strict preimage of $T_3$ produces
a section.
\bigbreak

In the case $s=-1$, the intersection numbers of the
line $T$ with the curve $R_3+T_3$ are always even;
then the preimage of $L$ is reducible and produces
two sections. We note also that $Q$ is in this
case an inflection point for both $R_3$ and $T_3$;
the common tangent line has also even intersection
numbers with $R_3+T_3$ and then it produces two sections.
We have found five different sections, then all of them.

\smallbreak
Let us consider now the second case. We know already
a section. By Proposition \copy\cover,
any other section should come from a section
to the pencil of lines through $S$ having always
even intersection numbers with the
ramification curve $R_3+T_3$. Then the problem is
as follows:

{\sl\narrower
Is there a curve $D$ of degree $d$ having a point of multiplicity
$d-1$ at $S$ and such that 
$(S\cdot R_3)_P\equiv(S\cdot T_3)_P \mod 2$ for any $P\in\PP^2$
and any branch of $D$ at $S$ has even intersection number with $T_3$?
\par}

Let us suppose that such a curve exists. It gives two different sections
$D_0$ and $D_1$ in the elliptic surface. From \cite{MP3},
$D_0$ and $D_1$ are torsion sections, and then they must be disjoint.
In particular, $D$ cannot intersect $R_3\cup T_3$ outside
$S,\AA_1,Q,Q'$ and no branch of $D$ at $S$ is tangent to
any branch of $T_3$ at $S$.
\bigbreak

$D_0$ and $D_1$ belong to the $5$-torsion, so by the
structure of the singular fibers, we have:

\smallbreak\item{--} $\AA_1\notin D$;

\smallbreak\item{--} $(T_3\cdot D)_{Q'}=(R_3\cdot D)_{Q'}=a=0,2,4$;

\smallbreak\item{--} $(T_3\cdot D)_{Q'}=(R_3\cdot D)_Q=b=1,3,5$.

\smallbreak
Then, putting all these conditions together, we obtain that
$S\notin D$ and so $D$ is a line; then $3=a+b$. The two
possibilities appear in the previous case, but not
in this one.
\qed
\enddemo

\subhead 3.- Case $m=11$
\endsubhead
\bigbreak

The method to find or discard the fibrations in the other
cases is the same one. As the answers are positive,
we will give the results that may be verified
by the reader. Let us consider the polynomial
$$
\multline
p_1(x,y,z):=\left ({\frac {11593}{95004009}}-
{\frac {4027\,v}{190008018}}\right ){
y}^{4}{x}^{2}+\left ({\frac {4705}{10556001}}-{\frac {2183\,v}{
10556001}}\right )zx{y}^{4}+\\
\left (-{\frac {1493\,v}{4691556}}+{\frac 
{803}{2345778}}\right ){z}^{2}{y}^{4}+\left (-{\frac {48226}{5000211}}
+{\frac {1475\,v}{5000211}}\right )z{y}^{3}{x}^{2}+\\
\left ({\frac {1174
\,v}{185193}}-{\frac {4736}{185193}}\right ){z}^{2}x{y}^{3}+\left ({
\frac {635\,v}{123462}}-{\frac {755}{61731}}\right ){z}^{3}{y}^{3}+\\
\left ({\frac {20153}{87723}}+{\frac {1081\,v}{175446}}\right ){z}^{2}
{y}^{2}{x}^{2}+
\left ({\frac {854}{3249}}-{\frac {187\,v}{3249}}
\right ){z}^{3}{y}^{2}x+\left (-{\frac {427}{6498}}+{\frac {187\,v}{
12996}}\right ){z}^{4}{y}^{2}+\\
\left (-{\frac {22612}{13851}}+{\frac {
386\,v}{13851}}\right ){z}^{3}y{x}^{2}+
\left ({\frac {1412}{1539}}+{
\frac {20\,v}{1539}}\right ){z}^{4}xy+{x}^{3}{z}^{3}+\left (-{\frac {
11\,v}{729}}-{\frac {485}{729}}\right ){z}^{4}{x}^{2}
\endmultline
$$
where ${v}^{2}+2=0$.

\proclaim{Proposition\cita} The curve $p_1(x,y,z)=0$ is an irreducible curve
with singularities $\EE_6$ (at $[1:0:0]$ and tangent line $z=0$), 
$\AA_1$ (at  $[0:0:1]$), $\AA_9$ (at $[0:1:0]$) and $\AA_2$ (at $[1:1:1]$).
The pencil of lines through the triple point determine after a double covering
an elliptic $K3$ fibration of type $[1,1,1,2,3,16]$ with trivial Mordell-Weil
group.
\endproclaim

\demo{Proof} The computations have been performed with MAPLEV. We note
that the curve is irreducible as the line $x=0$ joining $\AA_9$ and $\AA_1$
is not a component. Miranda-Persson classification finishes the  result.
\qed 
\enddemo

\subhead 4.- Case $m=13$
\endsubhead
\bigbreak

\proclaim{Proposition\cita} The curve $p_2(x,y,z)=0$ (see
below) is an irreducible curve
with singularities $\EE_6$ (at $[1:0:0]$ and tangent line $y=0$), 
$\AA_7$ (at  $[0:0:1]$), $\AA_4$ (at $[0:1:0]$) and $\AA_1$ (at $[1:1:1]$).
The pencil of lines through the triple point determine after a double covering
an elliptic $K3$ fibration of type $[1,1,1,2,5,14]$ with trivial Mordell-Weil
group.
\endproclaim

\demo{Proof} As before, computations have been performed with MAPLEV. We note
that the curve is irreducible as the line $x=y$ joining $\AA_7$ and $\AA_1$
is not a component. Miranda-Persson classification finishes the  result.
\qed 
\enddemo

We have:
$$
\multline
p_2(x,y,z):={y}^{3}{x}^{3}+\left (-{\frac {24284}{130321}}+{\frac
{10287\,v}{260642}}+{
\frac {144295\,{v}^{2}}{1824494}}\right ){y}^{4}{x}^{2}+\\
\left (-{\frac {
6071515\,{v}^{2}}{130321}}-{\frac {2851308\,v}{130321}}+{\frac {13668817}{
130321}}\right )z{x}^{2}{y}^{3}\\
+\left ({\frac {38660279\,v}{260642}}+{\frac {
161684215\,{v}^{2}}{521284}}-{\frac {179634441}{260642}}\right ){z}^{2}{x}^{2
}{y}^{2}+\\
\left (-{\frac {252208635\,{v}^{2}}{521284}}-{\frac {60782001\,v}{
260642}}+{\frac {277127879}{260642}}\right ){z}^{3}{x}^{2}y\\
+\left ({\frac {
55758423\,v}{521284}}+{\frac {460287135\,{v}^{2}}{2085136}}-{\frac {125694751
}{260642}}\right ){z}^{4}{x}^{2}+\\
\left (-{\frac {10473}{6859}}+{\frac {2326\,
v}{6859}}+{\frac {32860\,{v}^{2}}{48013}}\right )zx{y}^{4}+\\
\left (-{\frac {
361050\,{v}^{2}}{6859}}-{\frac {176895\,v}{6859}}+{\frac {1579285}{13718}}
\right ){z}^{2}x{y}^{3}+\\
\left ({\frac {725753\,v}{13718}}+{\frac {1458065\,{v
}^{2}}{13718}}-{\frac {1564472}{6859}}\right ){z}^{3}x{y}^{2}\\
+\left ({\frac {
1625477}{13718}}-{\frac {191737\,v}{6859}}-{\frac {3045105\,{v}^{2}}{54872}}
\right ){z}^{4}xy\\
+\left (-{\frac {268}{361}}+{\frac {141\,v}{722}}+{\frac {
3495\,{v}^{2}}{10108}}\right ){z}^{2}{y}^{4}+
\left ({\frac {825}{722}}-{
\frac {255\,v}{361}}-{\frac {1175\,{v}^{2}}{1444}}\right ){z}^{3}{y}^{3}\\
+\left (-{\frac {686}{361}}+{\frac {1099\,v}{1444}}+{\frac {6055\,{v}^{2}}{
5776}}\right ){z}^{4}{y}^{2},
\endmultline
$$
where $5\,{v}^{3}-4\,{v}^{2}-14\,v+14=0$.

\subhead 5.- Case $m=27$
\endsubhead
\bigbreak

In this cases we only state the result and give
the equation of the polynomial as the proofs are
very similar to the previous ones.

\bigbreak

\proclaim{Proposition\cita} The curve $p_3(x,y,z)=0$ (see below)
is an irreducible curve
with singularities $\EE_6$ (at $[0:0:1]$ and tangent line $y=0$), 
$\AA_3$ (at  $[1:0:0]$), $\AA_5$ (at $[0:1:0]$) and 
$\AA_4$ (at $[1:1:1]$).
The pencil of lines through the triple point determine after a double covering
an elliptic $K3$ fibration of type $[1,1,1,5,6,10]$ with trivial Mordell-Weil
group.
\endproclaim

We have
$$
\multline
p_3(x,y,z):=\left (-{\frac {200\,{v}^{2}}{297}}-{\frac {425}{297}}-{\frac
{110\,v}
{27}}\right ){y}^{4}{x}^{2}+\left ({\frac {125}{396}}+{\frac {5\,v}{9}
}-{\frac {13\,{v}^{2}}{396}}\right )z{y}^{4}x\\
+\left ({\frac {5\,{z}^{2
}}{528}}-{\frac {5}{264}}+{\frac {5\,v}{48}}\right ){z}^{2}{y}^{4}+
\left ({\frac {115\,{v}^{2}}{81}}+{\frac {220}{81}}+{\frac {875\,v}{81
}}\right ){y}^{3}{x}^{3}\\
+\left ({\frac {655}{108}}+{\frac {493\,v}{54}
}+{\frac {133\,{v}^{2}}{108}}\right )z{y}^{3}{x}^{2}+\left ({\frac {5
\,{v}^{2}}{36}}-{\frac {115}{36}}-{\frac {5\,v}{9}}\right ){z}^{2}{y}^
{3}x+{z}^{3}{y}^{3}\\
+\left (-{\frac {2225}{972}}-{\frac {3275\,v}{486}}
-{\frac {725\,{v}^{2}}{972}}\right ){y}^{2}{x}^{4}+\left (-{\frac {
2831}{324}}-{\frac {2032\,v}{81}}-{\frac {797\,{v}^{2}}{324}}\right )z
{y}^{2}{x}^{3}\\
+\left (-{\frac {37\,{v}^{2}}{72}}-{\frac {35}{36}}-{
\frac {215\,v}{72}}\right ){z}^{2}{y}^{2}{x}^{2}+\left ({\frac {1225\,
{z}^{2}}{972}}+{\frac {5215}{972}}+{\frac {7495\,v}{486}}\right )zy{x}
^{4}\\
+\left ({\frac {1105}{324}}+{\frac {788\,v}{81}}+{\frac {193\,{v}^
{2}}{324}}\right ){z}^{2}y{x}^{3}+\left (-{\frac {893\,{v}^{2}}{3888}}
-{\frac {4333}{1944}}-{\frac {24499\,v}{3888}}\right ){z}^{2}{x}^{4}
\endmultline
$$
where $25+75\,v+15\,{v}^{2}+{v}^{3}=0$.

\subhead 6.- Case $m=32$
\endsubhead

Let us consider the polynomial
$$
\multline
p_4(x,y,z):={y}^{3}{z}^{3}+\left ({\frac {5625\,v}{668168}}-{\frac
{33625}{334084}
}\right ){z}^{2}{x}^{4}+\left ({\frac {3475\,v}{58956}}+{\frac {39275}
{29478}}\right )y{z}^{2}{x}^{3}\\
+\left (-{\frac {1465\,v}{1734}}-{
\frac {1775}{867}}\right ){y}^{2}{x}^{2}{z}^{2}+\left ({\frac {173\,v}
{204}}-{\frac {299}{102}}\right ){y}^{3}x{z}^{2}+\left (-{\frac {v}{40
}}+{\frac {17}{20}}\right ){y}^{4}{z}^{2}\\
+\left ({\frac {19675\,v}{
501126}}-{\frac {188825}{501126}}\right )yz{x}^{4}+\left ({\frac {350
\,v}{4913}}+{\frac {23110}{4913}}\right ){y}^{2}{x}^{3}z+\left (-{
\frac {1580\,v}{867}}-{\frac {5900}{867}}\right ){y}^{3}{x}^{2}z\\
+
\left ({\frac {11\,v}{15}}-5/3\right ){y}^{4}xz+\left ({\frac {29555\,
v}{668168}}-{\frac {232705}{668168}}\right ){y}^{2}{x}^{4}+\left (-{
\frac {1885\,v}{29478}}+{\frac {116975}{29478}}\right ){y}^{3}{x}^{3}\\
+
\left (-{\frac {1205\,v}{1734}}-{\frac {33517}{8670}}\right ){y}^{4}{x
}^{2}
\endmultline
$$
where $v^2-v+34=0$.

\proclaim{Proposition\cita} The curve $p_4(x,y,z)=0$ is an irreducible curve
with singularities $\EE_6$ (at $[0:0:1]$ and tangent line $y=0$), 
$\AA_8$ (at  $[1:0:0]$), $\AA_2$ (at $[0:1:0]$) and two points
of type $\AA_1$ in the line $x+y+z=0$.
The pencil of lines through the triple point determine after a double covering
an elliptic $K3$ fibration of type $[1,1,2,2,3,15]$ with trivial Mordell-Weil
group.
\endproclaim

\subhead 7.- Case $m=37$
\endsubhead

\proclaim{Proposition\cita} The curve $p_5(x,y,z)=0$ (see below) 
is an irreducible curve
with singularities $\EE_6$ (at $[0:0:1]$ and tangent line $x=0$), 
$\AA_2$ (at  $[0:1:0]$), $\AA_8$ (at $[1:0:0]$) and two points
of type $\AA_1$ in the line $x+y+z=0$.
The pencil of lines through the triple point determine after a double covering
an elliptic $K3$ fibration of type $[1,1,2,2,9,9]$ with trivial Mordell-Weil
group.
\endproclaim

We have:
$$
\multline
p_5(x,y,z):=
\left ({\frac {3970803\,v}{130438}}-{\frac {345557847\,{v}^{2}}{65219}
}+{\frac {8058927}{130438}}\right ){y}^{4}{x}^{2}\\
+\left (-{\frac {
82574784\,{v}^{2}}{5929}}+{\frac {37159110\,v}{5929}}
-{\frac {3105297}
{5929}}\right )z{y}^{4}x\\
+\left (-{\frac {653967}{2156}}+{\frac {
3545235\,v}{1078}}-{\frac {5380479\,{v}^{2}}{1078}}\right ){z}^{2}{y}^
{4}\\
+\left ({\frac {5894214\,v}{9317}}
-{\frac {295704\,{v}^{2}}{9317}}-
{\frac {650011}{9317}}\right ){y}^{3}{x}^{3}\\
+\left (-{\frac {278076\,{
v}^{2}}{847}}+{\frac {808926\,v}{847}}-{\frac {86286}{847}}\right )z{y
}^{3}{x}^{2}\\
+\left (-{\frac {105723\,{v}^{2}}{77}}+{\frac {80505\,v}{
77}}-{\frac {15255}{154}}\right ){z}^{2}x{y}^{3}\\
+\left ({\frac {14286}
{1331}}-{\frac {136113\,v}{1331}}+{\frac {65742\,{v}^{2}}{1331}}
\right ){y}^{2}{x}^{4}+\left (-{\frac {24048\,v}{121}}+{\frac {30018\,
{v}^{2}}{121}}+{\frac {4599}{242}}\right )z{y}^{2}{x}^{3}\\
+\left (-{
\frac {2199\,v}{11}}+{\frac {3966\,{v}^{2}}{11}}+{\frac {195}{11}}
\right ){z}^{2}{y}^{2}{x}^{2}+\left (-{\frac {309}{121}}+{\frac {3711
\,v}{121}}-{\frac {8358\,{v}^{2}}{121}}\right )zy{x}^{4}\\
+\left ({
\frac {471\,v}{11}}-{\frac {903\,{v}^{2}}{11}}-{\frac {87}{22}}\right 
){z}^{2}y{x}^{3}+\left (-{\frac {42\,{v}^{2}}{11}}+{\frac {159\,v}{44}
}-{\frac {15}{44}}\right ){z}^{2}{x}^{4}+{z}^{3}{x}^{3}
\endmultline
$$
where $28\,{v}^{3}-30\,{v}^{2}+12\,v-1=0$.

\subhead 8.- Case $m=38$
\endsubhead
\bigbreak

Let us consider the polynomial
$$
\multline
p_6(x,y,z):=
{\frac {1404\,{x}^{2}{y}^{4}}{1445}}-{\frac {9\,x{y}^{4}z}{85}}+{
\frac {17\,{z}^{2}{y}^{4}}{60}}+{\frac {10800\,{x}^{3}{y}^{3}}{4913}}+
{\frac {1980\,{x}^{2}{y}^{3}z}{289}}\\
-{\frac {37\,{z}^{2}{y}^{3}x}{102}
}+{y}^{3}{z}^{3}
+{\frac {105840\,{x}^{4}{y}^{2}}{83521}}+{\frac {4410
\,{x}^{3}{y}^{2}z}{289}}+{\frac {13965\,{z}^{2}{y}^{2}{x}^{2}}{1156}}\\
+
{\frac {720300\,{x}^{4}yz}{83521}}+{\frac {780325\,{z}^{2}y{x}^{3}}{
29478}}
+{\frac {14706125\,{z}^{2}{x}^{4}}{1002252}}.
\endmultline
$$

\proclaim{Proposition\cita} The curve $p_6(x,y,z)=0$ is an irreducible curve
with singularities $\EE_6$ (at $[0:0:1]$ and tangent line $y=0$), 
$\AA_7$ (at  $[1:0:0]$), $\AA_1$ (at $[0:1:0]$) and two points
of type $\AA_2$ in the line $x+y+z=0$.
The pencil of lines through the triple point determine after a double covering
an elliptic $K3$ fibration of type $[1,1,2,3,3,14]$ with trivial Mordell-Weil
group.
\endproclaim

\subhead 9.- Case $m=55$
\endsubhead
\bigbreak

Let us consider the polynomial
$$
\multline
p_7(x,y,z):=
\left ({\frac {139}{176}}+{\frac {175\,v}{176}}\right ){y}^{4}{z}^{2}+.
\left (-{\frac {837\,v}{242}}+{\frac {7101}{968}}\right ){y}^{4}zx\\
+
\left ({\frac {30537}{10648}}-{\frac {29565\,v}{10648}}\right ){y}^{4}
{x}^{2}
+\left (-{\frac {151\,v}{44}}+{\frac {155}{44}}\right ){y}^{3}{
z}^{2}x+\left ({\frac {675}{242}}+{\frac {837\,v}{242}}\right ){y}^{3}
z{x}^{2}\\
+\left (-{\frac {669\,v}{2662}}+{\frac {2765}{1331}}\right ){y
}^{3}{x}^{3}+\left (-{\frac {81\,v}{22}}+{\frac {243}{44}}\right ){y}^
{2}{z}^{2}{x}^{2}+\left ({\frac {441\,v}{242}}-{\frac {183}{242}}
\right ){y}^{2}z{x}^{3}\\
+\left (-{\frac {1107}{1331}}+{\frac {2025\,v}{
1331}}\right ){y}^{2}{x}^{4}+\left (-{\frac {17}{11}}+{\frac {107\,v}{
22}}\right )y{z}^{2}{x}^{3}\\
+\left ({\frac {153\,v}{121}}+{\frac {18}{
121}}\right )yz{x}^{4}+{z}^{3}{x}^{3}+\left ({\frac {13}{22}}-{\frac {
5\,v}{22}}\right ){z}^{2}{x}^{4}
\endmultline
$$
where $3v^2-4v+2=0$.

\proclaim{Proposition\cita} The curve $p_7(x,y,z)=0$ is an irreducible curve
with singularities $\EE_6$ (at $[0:0:1]$ and tangent line $x=0$), 
$\AA_1$ (at  $[0:1:0]$), $\AA_7$ (at $[1:0:0]$) and two points
of type $\AA_2$ in the line $x+y+z=0$.
The pencil of lines through the triple point determine after a double covering
an elliptic $K3$ fibration of type $[1,1,3,3,8,8]$ with trivial Mordell-Weil
group.
\endproclaim

\finparrafo

\head\sec The complete determination of the Mordell-Weil group
for each type of semi-stable extremal fibrations
\endhead

In this section, we shall show Theorem \copy\thmi\
which will 
follow from the Table in \cite{MP3}, and the
Lemmas below.  We recall Lemma (1.3) and Shioda-Inose's
result that the isomorphism class of a $K3$ surface $X$ of Picard
number 20 is uniquely determined by the transcendental
lattice  $T_X$, modulo the action of  $SL_2({\bold Z})$
\cite{SI}.

\newbox\trans
\global\setbox\trans=\etiqueta
\proclaim{Lemma\cita} Let  $S$  be an even symmetric lattice of
rank $20$  and signature  $(1, 19)$  and
$T$  a positive definite even symmetric lattice of rank $2$.
Assume that
$\varphi : T^{\vee}/T \rightarrow S^{\vee}/S$  is an isomorphism
which induces the following equality involving
${\bold Q}/2{\bold Z}$-valued discriminant (quadratic) forms:
$$
q_{S} = - q_{T}.
$$
Let  $X$  be the unique $K3$ surface (up to isomorphisms) with
the transcendental lattice  $T_X = T$.  Then
the Picard lattice  $Pic X$  is isometric to  $S$.
\endproclaim

\demo{Proof} Consider the overlattice  $L$  of  $S \oplus T$  obtained
by adding all elements  $\varphi(x) + x$, $x \in T^{\vee}$,
where  $\varphi(x) \in S^{\vee}$  denotes one representative
of  $\varphi(x + T) \in S^{\vee}/S$.
The (even) intersection form on  $S \oplus T$  is
naturally extended to  a  ${\bold Q}$-valued one on
$S^{\vee} \oplus T^{\vee}$.  For each  $x \in T^{\vee}$, we have, modulo 
$2{\bold Z}$,
$(\varphi(x) + x, \varphi(x) + x) = (\varphi(x), \varphi(x)) + (x, x)
= q_{S}(\varphi(x)) + q_{T}(x) = 
-q_T(x) + q_T(x) = 0$,
i.e., $(\varphi(x) + x, \varphi(x) + x) \in 2{\bold Z}$.
Also for  $x_i \in T^{\vee}$,
combining  $(\varphi(x_1+x_2), \varphi(x_1+x_2)) = -(x_1 + x_2, x_1 + x_2)$
(mod $2{\bold Z}$) and  $(\varphi(x_i), \varphi(x_i)) = -(x_i, x_i)$ 
(mod $2{\bold Z}$), we see that  $(\varphi(x_1), \varphi(x_2)) = 
-(x_1, x_2)$ (mod ${\bold Z}$), whence  mod ${\bold Z}$
we have  $(\varphi(x_1) + x_1, \varphi(x_2) + x_2) = 
(\varphi(x_1), \varphi(x_2)) + (x_1, x_2) = 0$.
Thus  $L$  is an even (integral) symmetric 
lattice of rank 22 and signature
$(1+2, 19+0)$.  Clearly, $L/(S \oplus T) \cong T^{\vee}/T$
and hence  $|det(L)| = |det(S \oplus T)|/|T^{\vee}/T|^2 = 1$.
Now by the classification of indefinite unimodular 
even symmetric lattices,
$L$  is isometric to the $K3$ lattice (cf. \cite{Se}).

\smallbreak
On the other hand, by [SI], there is a unique $K3$
surface  $X$  (modulo isomorphisms) with the intersection form
of the transcendental lattice  $T_X$  equal to  $T$ (modulo
$SL_2({\bold Z})$).  We identify  $L$  with  $H^2(X, {\bold Z})$
and  $T$  with  $T_X$.  Note that there are two embeddings
$\iota_k : T_X \rightarrow H^2(X, {\bold Z})$:
$\iota_1 : T_X \hookrightarrow H^2(X, {\bold Z})$  as the
transcendental sublattice, and  $\iota_2 : T_X = T
\hookrightarrow S \oplus T \hookrightarrow L = H^2(X, {\bold Z})$.

\smallbreak
The embedding  $\iota_1$ (resp. $\iota_2$)  is primitive
by the definition of  $T_X$ (resp. of  $L$).
Now Nikulin's uniqueness theorem of primitive embedding
implies that there is an isometry  $\Psi$  of  $H^2(X, {\bold Z})$
such that  $\iota_1 = \Psi \circ \iota_2$
\cite{Mo, Cor.2.10}.  Note that the Picard lattice
$Pic X = (\iota_1(T_X))^{\perp} = (\Psi(\iota_2(T_X)))^{\perp}
= \Psi(T^{\perp}) = \Psi(S) \cong S$.
This proves the lemma.
\qed\enddemo

\newbox\caseiv
\global\setbox\caseiv=\etiqueta
\proclaim{Lemma\cita} Let  $f : X \rightarrow {\bold P}^1$
be of type  $m = 4$  as in Theorem \copy\thmi.  Then
$MW(f) \ne (0)$.
\endproclaim

\demo{Proof} 
Suppose the contrary that  $f : X \rightarrow {\bold P}^1$
is of type  $m = 4$  with  $MW(f) = (0)$.
Then  Pic $X$  is a direct sum  $U \oplus A_3 \oplus A_{15}$
of lattices, where  $U = (a_{ij})$  satisfies  
$a_{ii} = 0, a_{12} = a_{21} = 1$.  
Let  $(b_{ij})$  be the intersection matrix of
the transcendental lattice  $T = T_X$.
Then  $b_{ii} > 0$  and det$(b_{ij}) = |$det(Pic $X$)$| =
64$ (cf. \cite{BPV}). 
Modulo congruent action of  $SL(2, {\bold Z})$, we may
assume that  $-b_{11} < 2|b_{12}| \le b_{11} \le b_{22}$, 
and that $b_{12} \ge 0$  when  $b_{11} = b_{22}$.

\smallbreak
An easy calculation shows that one of the
following cases occurs:

\smallbreak
(1)  $(b_{ij}) =$ diag $[2, 32] \,\, $, (2)  $(b_{ij}) =$ diag $[4, 16],
\,\,$
\smallbreak 
(3)  $(b_{ij}) =$ diag $[8, 8], \,\,$  and
(4)  $b_{11} = 8, b_{22} = 10, b_{12} = 4$.

\smallbreak
Embed  $T$, as a sublattice, naturally  into  
$T^{\vee} = \text{\rm Hom}_{\bold Z}(T, {\bold Z})$.  
Then  $T^{\vee}/T \cong (\text{\rm Pic} X)^{\vee}/$ $(\text{\rm Pic} X)$ 
$\cong {\bold Z}/4{\bold Z} \oplus {\bold Z}/16{\bold Z}$.
Note that  $(\text{\rm Pic} X)^{\vee}/(\text{\rm Pic} X)$  is generated by
$\varepsilon_1 = (1/4)\sum_{i=1}^3 i v_i$  and  
$\varepsilon_2 = (1/16)\sum_{i=4}^{18} (i-3) v_i$, modulo
Pic $X$, where  $v_i$'s  form a canonical basis
of  $A_3 \oplus A_{15} \subseteq$ Pic $X$.
So the discriminantal quadratic form
$$q_T : T^{\vee}/T \rightarrow {\bold Q}/2{\bold Z}$$
is equal to  $-q_{\Pic X} = (-\varepsilon_1^2) \oplus 
(-\varepsilon_2^2) = (3/4) \oplus (15/16)$.

\smallbreak
On the other hand, in Case (4), $T^{\vee}$  has a  ${\bold Z}$-basis
$(e_1 \,\, e_2) (b_{ij})^{-1} = (g_1 \,\, g_2)$,
where  $e_1,e_2$  form a canonical basis of  $T$,
where  $g_1 = (1/32)(5e_1-2e_2), g_2 = (1/16)(-e_1+2e_2)$.
This leads to that  ord$(g_1)$  is equal to  32  in
$T^{\vee}/T$, a contradiction.

\smallbreak
In Cases (1)-(3)  where  $T =$ diag $[s, t]$,
with  $(s,t) = (2,32), (4,16)$  or  $(8,8)$, 
the discriminantal quadratic form  $q_T$  is equal to  $(1/s) \oplus (1/t)$.
This leads to that  $(1/s) \oplus (1/t) \cong (3/4) \oplus (15/16)$, 
which is impossible by an easy check. 

\smallbreak
Therefore, the lemma is true.
\qed\enddemo
\newbox\casev
\global\setbox\casev=\etiqueta
\proclaim{Lemma\cita}  Consider the pairs below:
$$(m, G_m) = (2, \langle 0 \rangle), (9, \langle 0 \rangle), (11, \langle 0
\rangle),
(13, \langle 0 \rangle), (27, \langle 0 \rangle), (32, \langle 0 \rangle),$$
$$(37, \langle 0 \rangle), (38, \langle 0 \rangle),
(55, \langle 0 \rangle), (35, {\bold Z}/2{\bold Z}), 
(53, \langle {\bold Z}/3{\bold Z} \rangle).$$

\smallbreak 
For each of these eleven pairs  $(m, G_m)$,
there is a Jacobian elliptic  $K3$  surface
$f_m : X_m \rightarrow {\bold P}^1$  of type  $m$  as in Theorem \copy\thmi\
such that  $(m, MW(f_m)) = (m, G_m)$.
\endproclaim

\demo{Proof}
The existence of the pairs where
$m = 2, 35$  is proved constructively in \cite{AT}.
The rest is also constructively proved in \S 2.
In the paragraphs below, 
we will give an independent lattice-theoretical proof.

\smallbreak
Let  $T_m$, $m = 2, 9, 11, 13, 27, 32, 37, 38, 55, 35, 53$,
be the positive definite even symmetric lattice of rank 2 with 
the following intersection form, respectively:
$$\pmatrix 4&2 \\
           2&10 \endpmatrix,
\pmatrix 10&0 \\
           0&10 \endpmatrix,
\pmatrix 10&2 \\
         2&10 \endpmatrix,
\pmatrix 2&0 \\
         0&70 \endpmatrix,
\pmatrix 10&0 \\
         0&30 \endpmatrix,
\pmatrix 12&6 \\
         6&18 \endpmatrix,$$
$$\pmatrix 18&0 \\
         0&18 \endpmatrix,
\pmatrix 6&0 \\
         0&42 \endpmatrix,
\pmatrix 24&0 \\
         0&24 \endpmatrix,
\pmatrix 6&0 \\
         0&12 \endpmatrix,
\pmatrix 4&0 \\
         0&12 \endpmatrix.$$

For the first nine  $m$  above,
let  $S_m$  be the even lattice of rank 20 and
signature (1,19) with the following intersection form, respectively
$$U \oplus A_1 \oplus A_{17},
U \oplus A_9 \oplus A_9,
U \oplus A_1 \oplus A_2 \oplus A_{15},$$
$$U \oplus A_1 \oplus A_4 \oplus A_{13},
U \oplus A_4 \oplus A_5 \oplus A_9,
U \oplus A_1 \oplus A_1 \oplus A_2 \oplus A_{14},$$
$$U \oplus A_1 \oplus A_1 \oplus A_8 \oplus A_8,
U \oplus A_1 \oplus A_2 \oplus A_2 \oplus A_{13},
U \oplus A_2 \oplus A_2 \oplus A_7 \oplus A_7.$$

We now define  $S_m$  for  $m = 35, 53$.
Let  $\Gamma_{35}$  be the lattice
$U \oplus A_1 \oplus A_1 \oplus A_5 \oplus A_{11}$,
with  $G, H, J_i (1 \le i \le 5), \theta_i (1 \le i \le 11)$
as the canonical basis of   $A_1 \oplus A_1 \oplus A_5 \oplus A_{11}$,
and  ${\Cal O}, F$  as a basis of  $U$  such that  ${\Cal O}^2 = -2,
F^2 = 0, {\Cal O}\cdot F = 1$.  

\smallbreak
We extend  $\Gamma_{35}$
to an index-2 integral over lattice  
$S_{35} = \Gamma_{35} + {\bold Z} s_{35}$, where
$$
s_{35} = {\Cal O} + 2F - G/2 - H/2 -
(1/2)(\sum_{i=1}^6 i \theta_i + \sum_{i=7}^{11} (12-i) \theta_i).
$$
It is easy to see that the intersection form on  $\Gamma_{35}$
can be extended to an integral even symmetric lattice of signature
$(1, 19)$.  Indeed, setting  $s = s_{35}$, we have
$$
s^2 = -2, s \cdot F = s \cdot G = s \cdot H = s \cdot \theta_6 = 1,
s \cdot {\Cal O} = s \cdot J_i = s \cdot \theta_j = 0 
\,\, (\forall i; j \ne 6).
$$

Moreover, $|\det(S_{35})| = |\det(\Gamma_{35})|/2^2 = 72$.

\smallbreak
Note that  $\Gamma_{35}^{\vee} = \Hom_{\bold Z} (\Gamma_{35}, {\bold Z})$
contains naturally  $\Gamma_{35}$  as a sublattice  
with  ${\bold Z}/2{\bold Z} \oplus 
{\bold Z}/2{\bold Z} \oplus {\bold Z}/6{\bold Z} \oplus 
{\bold Z}/12{\bold Z}$  as the factor group, and 
is generated by the following, modulo  $\Gamma_{35}$:
$$h_1 = G/2, \,\, h_2 = H/2, \,\, h_3 = (1/6)\sum_{i=1}^5 i J_i, \,\,
h_4 = (1/12)\sum_{i=1}^{11} i \theta_i.$$

\smallbreak
Since  $(S_{35})^{\vee}$  is an (index-2) sublattice of  
$(\Gamma_{35})^{\vee}$, an element  $x$  is in  $(S_{35})^{\vee}$
if and only if  $x = \sum_{i=1}^4 a_i h_i$ (mod $\Gamma_{35}$)
such that  $x$  is integral on  $S_{35}$, i.e.,
$x \cdot s = (a_1 + a_2 + a_4)/2$  is an integer.
Hence  $(S_{35})^{\vee}$  is generated by the following, modulo  $\Gamma_{35}$:
$$
h_3, \,\, 2h_i, \,\, h_1 + h_2, \,\, h_1 + h_4, \,\, h_2 + h_4.
$$

\smallbreak
Noting that  $2h_1, 2h_2 \in S_{35}$  and
$(h_1+h_2) + 6h_4$  is equal to  $s$  (mod  $\Gamma_{35}$)
and hence contained in  $S_{35}$, we can see easily that
$(S_{35})^{\vee}$  is generated by the following, modulo  $S_{35}$:
$$
\varepsilon_1 := h_3, \,\,\,\, \varepsilon_2 := h_1-h_4.
$$
Now the fact that  $|(S_{35})^{\vee}/S_{35}| = 72$
and that  $6\varepsilon_1, 12\varepsilon_2 \in S_{35}$  imply that
$(S_{35})^{\vee}/S_{35}$  is a direct sum of its cyclic subgroups
which are of order 6, 12, and generated by
$\varepsilon_1, \varepsilon_2$, modulo  $S_{35}$.

\smallbreak
We note that the negative of the discriminant form
$$
-q_{(S_{35})} = (-(\varepsilon_1)^2) \oplus (-(\varepsilon_2)^2)
= (5/6) \oplus ((1/2) + (11/12)) = (5/6) \oplus (-7/12).
$$

\smallbreak 
Next we define  $S_{53}$.
Let  $\Gamma_{53}$  be the lattice
$U \oplus A_2 \oplus A_2 \oplus A_3 \oplus A_{11}$,
with  $G_i (i = 1,2), H_i (i = 1,2), J_i (i = 1, 2, 3), 
\theta_i (1 \le i \le 11)$
as the canonical basis of   $A_2 \oplus A_2 \oplus A_3 \oplus A_{11}$,
and  ${\Cal O}, F$  as a basis of  $U$  as in the case of  $S_{35}$.

\smallbreak
Extend  $\Gamma_{53}$  to an index-3 integral over lattice
$S_{53} = \Gamma_{53} + {\bold Z} s_{53}$, where
$$s_{53} = {\Cal O} + 2F - (1/3)(2G_1 + G_2 + 2H_1 + H_2) -
(2/3)\sum_{i=1}^{11} i \theta_i - \sum_{i=5}^{11} (i-4) \theta_i.$$
The intersection form on  $\Gamma_{53}$
can be extended to an integral even symmetric lattice of signature
$(1, 19)$  such that the following is true, where we set  $s = s_{53}$:
$$
s^2 = -2, s \cdot F = s \cdot G_1 = s \cdot H_1 = s \cdot \theta_4 = 1,
$$
$$
s \cdot {\Cal O} = s \cdot G_2 = s \cdot H_2 
= s \cdot J_i = s \cdot \theta_j = 0 \,\, 
(\forall i; j \ne 4).
$$

Moreover, $|\det(S_{53})| = |\det(\Gamma_{53})|/3^2 = 48$.

\smallbreak
Note that  $\Gamma_{53}^{\vee}$
contains naturally  $\Gamma_{53}$  as a sublattice
with  ${\bold Z}/3{\bold Z} \oplus
{\bold Z}/3{\bold Z} \oplus {\bold Z}/4{\bold Z} \oplus
{\bold Z}/12{\bold Z}$  as the factor group, and
is generated by the following, modulo  $\Gamma_{53}$:
$$h_1 = (1/3)\sum_{i=1}^2 i G_i, \,\,
h_2 = (1/3)\sum_{i=1}^2 i H_i, \,\,
h_3 = (1/4)\sum_{i=1}^3 i J_i, \,\,
h_4 = (1/12)\sum_{i=1}^{11} i \theta_i.$$

\par
Since  $(S_{53})^{\vee}$  is an (index-3) sublattice of
$(\Gamma_{53})^{\vee}$, an element  $x$  is in  $(S_{53})^{\vee}$
if and only if  $x = \sum_{i=1}^4 a_i h_i$ (mod $\Gamma_{53}$)
such that  $x$  is integral on  $S_{53}$, i.e.,
$x . s = (a_1 + a_2 + a_4)/3$  is an integer.
Hence  $(S_{53})^{\vee}$  is generated by the following, modulo  $\Gamma_{53}$:
$$h_3, \,\, 3h_i, \,\, h_1+h_2+h_4, \,\, h_1-h_2, \,\, h_1-h_4, \,\, h_2-h_4.$$

\par
Noting that  $3h_1, 3h_2 \in S_{53}$  and
$3h_4 + (h_1+h_2+h_4)$  is equal to  $s$  (mod  $\Gamma_{53}$)
and hence contained in  $S_{53}$, we see that
$(S_{53})^{\vee}$  is generated by   
$\varepsilon_1 := h_3, \varepsilon_2 := h_1-h_4$, modulo  $S_{53}$.
As in the case of  $S_{35}$,
$(S_{53})^{\vee}/S_{53}$  is a direct sum of its cyclic subgroups,
which are of order  4, 12, and generated by
$\varepsilon_1, \varepsilon_2$, modulo  $S_{53}$.

\par
The negative of the discriminant form
$$-q_{(S_{53})} = (-(\varepsilon_1)^2) \oplus (-(\varepsilon_2)^2)
= (3/4) \oplus ((2/3) + (11/12)) = (3/4) \oplus (-5/12).$$
\newbox\zerosec
\global\setbox\zerosec=\etiqueta
\proclaim{Claim\cita} The pair  $(S_m, T_m)$  satisfies the
conditions of Lemma (3.1) and hence if we let 
$X_m$  be the unique  $K3$  surface with  $T_{X_m} = T_m$
then  $\Pic X_m = S_m$ (both two equalities here are modulo
isometries).
\endproclaim

\demo{Proof of the claim}
We need to show that  $q_{T_m} = -q_{S_m}$.  Note that  
$A_n^{\vee}/A_n = {\bold Z}/(n+1){\bold Z}$
and  $q_{(A_n)} = (-n/(n+1))$. For the first nine $m$,
if we write  $S_m = U \oplus A_{n_1-1} \oplus \cdots A_{n_k-1}$,
then
$$
q_{S_m} = (-(n_1-1)/n_1) \oplus \cdots \oplus (-(n_k-1)/n_k);
$$
moreover, $S_m^{\vee}/S_m$  is generated 
by two elements  $\varepsilon_i$ ($i = 1, 2$) 
($\varepsilon_i$  is a simple sum of the natural generators of 
$S_m^{\vee}/S_m$)
such that for every  $a, b \in {\bold Z}$
one has  $-q_{(S_m)}(a\varepsilon_1 + a\varepsilon_2) = 
-a^2(\varepsilon_1)^2 - b^2(\varepsilon_2^2)$.
For all eleven  $m$, $\varepsilon_i$  can be chosen such that
$(-\varepsilon_1^2, -\varepsilon_2^2)$  
is respectively given as follows:
$$(1/2, 17/18), (9/10, 9/10), (1/2, -19/48), (1/2, 121/70),$$
$$(9/10, 49/30), (-5/6, -17/30),(25/18, 25/18), (-5/6, -17/42),$$
$$(-11/24, -11/24), (5/6, -7/12), (3/4, -5/12).$$

\smallbreak
On the other hand, $T_m^{\vee}$  is generated by
$(g_1 \,\, g_2) = (e_1 \,\, e_2) T_m^{-1}$, where
$e_1, e_2$  form a canonical basis of  $T_m$  which gives
rise to the intersection matrix of  $T_m$  shown before this claim.
Now, the claim follows from the existence of the following isomorphism,
which induces  $q_{T_m} = -q_{S_m}$: 
$$
\varphi : T_m^{\vee}/T_m \rightarrow S_m^{\vee}/S_m
$$
$$
(g_1 \,\, g_2) \mapsto (\varepsilon_1 \,\, \varepsilon_2) B_m.
$$
Here  $B_m$  is respectively given as:
$$\pmatrix 1&1 \\ 2&5 \endpmatrix,
\pmatrix 7&0 \\ 0&7 \endpmatrix,
\pmatrix 0&1 \\ 11&17 \endpmatrix,
\pmatrix 1&0 \\ 0&51 \endpmatrix,
\pmatrix 7&0 \\ 0&17 \endpmatrix,
\pmatrix -2&1 \\ 1&3 \endpmatrix,$$
$$\pmatrix 7&0 \\ 0&7 \endpmatrix,
\pmatrix 2&3 \\ 21&10 \endpmatrix,
\pmatrix 2&3 \\ 3&-2 \endpmatrix,
\pmatrix 3&2 \\ 4&3 \endpmatrix,
\pmatrix 0&1 \\ 3&4 \endpmatrix.$$
\enddemo

Write  $S_m$ (resp. $\Gamma_m$) as
$U \oplus {\bold A}(m)$  with
${\bold A}(m) = A_{n_1-1} \oplus \cdots \oplus A_{n_k-1}$, for the first 
nine  $m$ (resp. $m = 35, 53$) as in the definitions of them.
Let  ${\Cal O}, F$  be a  ${\bold Z}$-basis of  $U$
for all  $m$, as in the definition of  $S_{35}$.
By \cite{PSS, p. 573, Th 1}, after an (isometric) action of reflections
on  $S_m = \Pic X_m$, we may assume at the beginning 
that  $F$  is a fiber of an elliptic fibration  
$f_m : X_m \rightarrow {\bold P}^1$.  
Since  ${\Cal O}^2 = -2$, Riemann-Roch Theorem implies that
${\Cal O}$  is an effective divisor because  ${\Cal O} \cdot F > 0$.
Moreover, ${\Cal O}\cdot F = 1$  implies that 
${\Cal O} = {\Cal O}_1 + F'$  where  ${\Cal O}_1$  is 
a cross-section of  $f_m$  and  $F'$  is an effective divisor 
contained in fibers.  So  $f_m$  is a Jacobian elliptic
fibration and we can choose  ${\Cal O}_1$  as the zero element of
$MW(f_m)$.

\smallbreak
Let  $\Lambda_m$  be the lattice generated by 
all fiber components of  $f_m$.  Clearly, 
$\Lambda_m = {\bold Z} F \oplus \Delta,
\Delta = \Delta(1) \oplus \cdots
\oplus \Delta(r)$ (depending on  $m$), where
each  $\Delta(i)$  is a negative definite even lattice
of Dynkin type  $A_p, D_q$, or  $E_r$, contained in a single 
reducible singular fiber  $F_i$  of  $f_m$  and spanned by smooth
components of  $F_i$  disjoint from  ${\Cal O}_1$.

\newbox\onesec
\global\setbox\onesec=\etiqueta
\proclaim{Claim\cita} We have:

\smallbreak\item{\rm(1)} $Span_{\bold Z} \{x \in S_m |
x \cdot F = 0, x^2 = -2\}$ $= \Lambda_m$
$= {\bold Z} F \oplus {\bold A}(m)$;
in particular, $r = k$, and there are
lattice-isometries: $\Delta \cong {\bold A}(m)$
and  $\Delta(i) \cong A_{n_i}$  ($i = 1, 2, \dots, k$),
after relabelling.

\smallbreak\item{\rm(2)} There are  $k$  singular fibers  $F_i$  of type  
${\widetilde A}_{n_i-1}$ ($1 \le i \le k$)  of  $f_m$, 
and any fiber other than  $F_i$  is irreducible.

\smallbreak\item{\rm(3)} $MW(f_m) = (0)$ (resp. ${\bold Z}/2{\bold Z}$, ${\bold
Z}/3{\bold Z}$)
for the first nine  $m$ (resp. $m = 35, 53$).
\endproclaim

\demo{Proof}
The assertion (2) follows from (1) (see also
\cite{K, Lemma 2.2}).

The first equality in  (1) is clear from Kodaira's classification of
elliptic fibers and the Riemann Roch Theorem as used prior to
this claim to deduce  ${\Cal O} \ge 0$.
The second equality is clear for the cases of the first 
nine  $m$  because then  $\Pic X_m = S_m = ({\bold Z}{\Cal O} + 
{\bold Z}F) \oplus {\bold A}(m)$.
\smallbreak

Let  $m = 35, 53$.  We now show the second
equality using Lemma \copy\lemmasi.  
Clearly, ${\bold Z} F \oplus {\bold A}(m)$  is contained 
in the first term of (1) and hence in  $\Lambda_m$.  
One notes that  $19 =$ rank $S_m - 1 \ge$ rank $\Lambda_m =
1 +$ rank $\Delta \ge$ $1 +$ rank ${\bold A}(m) =$ 
$1 + \sum_{i=1}^k (n_i - 1) = 19$.  Hence  
$\Delta = \Delta(1) \oplus \cdots \oplus \Delta(r) \cong 
\Lambda_m/{\bold Z} F$  contains a finite-index sublattice
$({\bold Z} F \oplus {\bold A}(m)) /{\bold Z} F \cong {\bold A}(m) =$ 
$A_{n_1-1} \oplus \cdots \oplus A_{n_k-1}$.

\smallbreak
Suppose the contrary that the second equality in
(1) is not true.  Then  ${\bold A}(m)$  is an index-$n$
($n > 1$) sublattice of  $\Delta$.
By Lemma \copy\lemmasi, one of Cases (2-1) - (2-3) there occurs.

\smallbreak
Case (2-1).  Then  $m = 35$, $f_m$  has reducible singular fibers
of types  ${\widetilde A}_1, I_{13}^*$
and no other reducible fibers.  This leads to
that  $72 = |\Pic X_m| = (2 \times 4)/|MW(f_m)|$,
a contradiction (cf. \cite{S}).

\smallbreak
Case (2-2).  Then  $m = 53$, $f_m$  has reducible singular fibers
of types  ${\widetilde A}_2, I_{12}^*$
and no other reducible fibers.  This leads to
that  $48 = |\Pic X_m| = (3 \times 4)/|MW(f_m)|$,
a contradiction.

\smallbreak
Case (2-3).  Then  $m = 35$, $f_m$  has reducible singular fibers
of types  ${\widetilde A}_1, I_{12}, IV^*$
and no other reducible fibers.  Since 
$72 = |\Pic X_m| = (2 \times 12 \times 3)/|MW(f_m)|$,
we have  $MW(f_m) = (0)$  and  $S_m = \Pic X_m =$ {\bf Z}${\Cal O}_1 
+ \Lambda_m =$ {\bf Z}${\Cal O}_1 + ($ {\bf Z}$F \oplus \Delta) =$ 
{\bf Z}${\Cal O}_1 +($ {\bf Z}$F \oplus A_1 \oplus A_{11} \oplus E_6)$.

By the Riemann-Roch theorem and the fact that
$(s_m)^2 = -2$, $s_m . F = 1$  and  $MW(f_m) = (0)$,
we see that  $s_m = {\Cal O}_1$ (mod $\Lambda_m$).  
This, together with
the fact that  ${\Cal O} = {\Cal O}_1$ (mod $\Lambda_m$)
and the definition of  $s_m$, implies that
$(1/2)(G + H + D) \in \Lambda_m$, where  $D = 
\sum_{i=1}^6 i \theta_i + \sum_{i=7}^{11} (12-i) \theta_i$.

Consider the index-2 extension  
$$
A_1 \oplus A_{11} \oplus
(A_1 \oplus A_5)\! = \!{\bold A}(m)\! \cong \!
({\bold Z}F \oplus {\bold A}(m))/{\bold Z}F
\!\subseteq\! ({\bold Z}F \oplus \Delta)/{\bold Z}F\! \cong\! \Delta \!= 
\! A_1 \oplus A_{11} \oplus E_6.
$$  
The proof of Lemma \copy\lemmasi\ shows that  
(the first summand  $A_1$  in this rearranged  ${\bold A}(m)$)
$\oplus {\bold Z}F =$ (the summand  $A_1$  
in  $\Delta$) $\oplus {\bold Z}F$,
(the summand  $A_{11}$  in  ${\bold A}(m)$)
$\oplus {\bold Z}F =$ (the summand  $A_{11}$  in  $\Delta$)
$\oplus {\bold Z}F$, and (the summand  $(A_1 \oplus A_5)$  in  ${\bold A}(m)$)
$\oplus {\bold Z}F \subseteq$ (the summand  $E_6$  in  $\Delta$)
$\oplus {\bold Z}F$.  So we may assume that, mod ${\bold Z} F$,
$G$  is the {\bf Z}-generator of the first summand  $A_1$ 
in  $\Delta$, $\theta_i$ ($1 \le i \le 11$)
form a  {\bf Z}-basis of the summand  $A_{11}$  in  $\Delta$,  
and  $H$  is contained in the summand  $E_6$  in  $\Delta$. 

In particular, for  $(G + H + D)/2 \in \Lambda_m = {\bold Z}F \oplus \Delta =
{\bold Z}F \oplus (A_1 \oplus A_{11} \oplus E_6)$,
we have, mod ${\bold Z} F$, $G/2 \in A_1$, 
$H/2 \in E_6$, and  $D/2 \in A_{11}$.
We reach a contradiction to the above observation 
that the  $A_1$  in  $\Delta$  is generated by
$G$  over  {\bf Z}.
 
\smallbreak
Therefore, the second equality of (1) is true.  So
there is an isometry  $\Phi : \Delta \cong \Lambda_m/{\bold Z} F
\cong {\bold A}(m)$.  Now the rest of (1) follows from Lemma \copy\lemmasi.

The assertion (3) follows from the fact in 
\cite{S, Th 1.3},
that  $MW(f_m)$  is isomorphic to
the factor group of  $\Pic X_m$  modulo
$({\bold Z}{\Cal O}_1 + {\bold Z}F) \oplus \Delta$, 
where the latter is equal to
$({\bold Z}{\Cal O} + {\bold Z}F) + \Delta =
({\bold Z}{\Cal O} + {\bold Z}F) \oplus {\bold A}(m) =$
$U \oplus {\bold A}(m)$.
This proves the claim.
\enddemo

\smallbreak
The existence of singular fibers  $F_i$ ($i = 1,2, \dots, k$)
of type  $I_{n_i-1}$,
the fact that the sum of Euler numbers of singular fibers
of  $f_m$  is 24, the fact that each fiber other than  $F_i$
is irreducible, and \cite{MP3, Lemma 3.1 and Proposition 3.4} 
imply that  $f_m$  is semi-stable.  Hence  $F_i$ ($i = 1,2, \dots, k$)  
is of type  $I_{n_i}$, there are 
$\chi(X_m) - \sum_i (n_i - 1) - k = 6 - k$
fibers of type  $I_1$, and  $f_m$  is of type
$[1, 1, \dots, 1, n_1, \dots, n_k]$, i.e., of type  $m$
after an easy case-by-case check.
Moreover, $(m, MW(f_m)) = (m, G_m)$
for all eleven  $m$  by the last claim.
This completes the lattice-theoretical proof of Lemma \copy\casev.
\qed\enddemo

\smallbreak
\remark{Remark\cita} We note that  $S_{35}
= U \oplus A_1 \oplus A_{11} \oplus E_6$.  This is because
the lattices  $T_{35}$  and the one on the right hand side 
satisfy all conditions of Lemma \copy\trans\ by an easy check.
In particular, using \cite{MP3, Lemma 3.1 and Proposition 3.4}
as in the proof of Lemma \copy\casev, we can show that
there is a Jacobian elliptic fibration
$\tau_m : X_m \rightarrow {\bold P}^1$ ($m = 35$)
with singular fibers  $I_1, I_1, I_2, I_{12}, IV^*$
and with  $MW(\tau_m) = (0)$.
\endremark

\newbox\lemmaii
\global\setbox\lemmaii=\etiqueta
\proclaim{Lemma\cita} Let  $f : X \rightarrow {\bold P}^1$
be of type  $m$  as in Theorem \copy\thmi.  Then the following are true:

\smallbreak\item{\rm(1)} If  $m = 48$, then  $MW(f) \ne {\bold Z}/2{\bold Z}$, 
or  ${\bold Z}/4{\bold Z}$.

\smallbreak\item{\rm(2)} If  $m = 4$, then  $MW(f) \ne {\bold Z}/2{\bold Z}$.

\smallbreak\item{\rm(3)} If  $m = 31$, then  $MW(f) \ne {\bold Z}/2{\bold Z}$.

\smallbreak\item{\rm(4)} If  $m = 44$, then  $MW(f) \ne {\bold Z}/2{\bold Z}$.

\smallbreak\item{\rm(5)} If  $m = 69$, then it is impossible that  
$MW(f)$  is  ${\bold Z}/2{\bold Z}$  with  $s = (0,0,0,0,2,6)$
as its generator (see Remark \copy\gnrt).

\smallbreak\item{\rm(6)} If  $m = 92$, then  $MW(f) \ne {\bold Z}/2{\bold Z}$.
\endproclaim

\demo{Proof} Let  $f : X \rightarrow {\bold P}^1$  be
of type  $m$  as in Theorem \copy\thmi.

\smallbreak
(1)  Assume that  $f$  is of type  $m = 48$ 
and  $MW(f) \supseteq {\bold Z}/2{\bold Z}$.  
We will show that  $MW(f) \supseteq {\bold Z}/8{\bold Z}$
which will imply (1).

\smallbreak
$m = 48$  means that the singular fiber type of  $f$  is
$I_1, I_1, I_2, I_4, I_8, I_8$.
Using the height pairing in \cite{S} or the Table in \cite{MP3}, 
we may assume that
$MW(f)$  contains  $s = (0,0,0,0,4,4)$  as a 2-torsion section
after suitable labeling of fibre components.

\smallbreak
Let  $Y$, a $K3$ surface again, be the minimal resolution 
of the quotient surface  $X/\langle s \rangle$.  
$f$  on  $X$  induces a Jacobian semi-stable elliptic fibration
$g : Y \rightarrow {\bold P}^1$  of singular fiber type
$I_2, I_2, I_4, I_8, I_4, I_4$  where these 6 ordered
singular fibers are respectively ``images'' of ordered
singular fibers on  $X$.

\smallbreak
To be precise, let  $\sigma : {\widetilde X} \rightarrow X$
be the blowing-up of all 8 intersections in the
first 4 singular fibers of  $f$  of types  $I_1, I_1, I_2, I_4$.
Then  $Y = {\widetilde X}/\langle s \rangle$
and the  ${\bold Z}/2{\bold Z}$-covering  
$\pi : {\widetilde X} \rightarrow Y$  is branched
along 4 disjoint curves  $\theta^{(i)}_j$,
where  $(i,j) = (1,1), (2,1), (3,1), (3,3), (4,1),
(4,3)$, $(4,5), (4,7)$.  Here we choose the common
image of the zero section and the 2-torsion section
$s$  of  $f$, as the zero section  $O_1$  of  $g$,
and label clock or anti-clock wise the  $i$-th singular 
fiber of  $g$  of type  $I_{n_i}$  as  $\sum_{j=0}^{n_i-1}
\theta^{(i)}_j$  so that  $O_1$  passes through  $\theta^{(i)}_0$,
where  $[n_1, \dots, n_6] = [2,2,4,8,4,4]$.

\smallbreak
Note that  $(Y, g)$  is of type $m = 103$  in the Table of \cite{MP3}
and hence there is a 4-torsion section  $t$  of  $g$
equal to  $(0,0,2,2,1,1)$  or  $(0,0,1,2,1,2)$  or  $(0,0,1,2,2,1)$,
after choosing either clockwise or counterclockwise 
labeling of fiber components,  
where for orders of six fibers of  $g$  we use the current indexing 
inheriting from that of  $f$.

\smallbreak
If   $t = (0,0,1,2,1,2)$  or  $(0,0,1,2,2,1)$, 
then  $t$  meets the branch
locus of  $\pi$  transversally at one point only so that 
$\pi^{-1}(t)$  is a smooth irreducible curve
and  $\pi : \pi^{-1}(t) \rightarrow t$  is a double cover
with exactly one ramification point, a contradiction
to Hurwitz's genus formula applied to the covering map  $\pi$.

\smallbreak
Thus  $t = (0,0,2,2,1,1)$.  A check using height pairing
shows that  $\pi^{-1}(t)$  is a disjoint union of
two 8-torsion sections of  $f$.  Hence  
$MW(f) \supseteq {\bold Z}/8{\bold Z}$.
Indeed, $MW(f) = {\bold Z}/8{\bold Z}$  by \cite{MP3}.
This proves (1).

\smallbreak 
Now assume that  $f$  is of type  $m = 4$
(resp. $m = 31$, $m = 44$, $m = 69$ 
with  $MW(f) = \langle s = (0,0,0,0,2,6) \rangle$), or  $m = 92$)
and  $MW(f) \supseteq {\bold Z}/2{\bold Z}$.
Then  $MW(f)$  contains a unique 2-torsion section  
$s = (0,0,0,0,0,8)$ (resp. $s = (0,0,0,0,0,8)$,
$s = (0,0,0,0,2,6)$, $s = (0,0,0,0,2,6)$, or  $s = (0,0,0,2,2,4)$)
(cf. \cite{MP3}). 
As in (1) we can show that  $f$  induces a
Jacobian semi-stable elliptic fibration  $g$
on the minimal resolution  $Y$  of  $X/\langle s \rangle$.
The singular fiber type of  $g$  is  $I_{n_1} + \dots + I_{n_6}$
where  $[n_1, \dots, n_6]$  is equal to
$[2,2,2,2,8,8]$ (resp. $[2,2,4,4,4,8]$, $[2,2,4,8,2,6]$,
$[2,4,4,6,2,6]$, or  $[2,6,8,2,2,4]$)
and hence is of type  $m = 94$ (resp. $m = 103$, $m = 97$,
$m = 104$, or  $m = 97$) in the Table of  
\cite{MP3}.  Now the inverse on  $X$  of the 2-torsion section
$t = (0,0,0,0,4,4)$  (resp. $t = (0,0,0,2,2,4)$, $t = (0,0,0,4,1,3)$,
$t$  is one of  $(0,2,2,0,1,3)$  and  $(1,2,2,3,0,0)$, 
or  $t = (0,0,0,4,1,2)$)  on  $Y$  is a disjoint union of
two 4-torsion sections of  $f$.  Hence  
$MW(f) \supseteq {\bold Z}/4{\bold Z}$.
Indeed, $MW(f) = {\bold Z}/4{\bold Z}$  by \cite{MP3}.
This proves (2) - (6).  The proof of the lemma is completed.
\enddemo

\proclaim{Lemma\cita} Let  $f : X \rightarrow {\bold P}^1$
be of type  $m$  as in Theorem \copy\thmi.  Then
each of the following pairs  $(m, MW(f))$  occurs:

\smallbreak
$$(69, {\bold Z}/2{\bold Z} = \langle (0,1,1,0,0,6) \rangle),
(69, {\bold Z}/4{\bold Z}),
(92, {\bold Z}/4{\bold Z}),$$
$$(32, {\bold Z}/3{\bold Z}), 
(37, {\bold Z}/3{\bold Z}), 
(44, {\bold Z}/4{\bold Z}), 
(55, {\bold Z}/2{\bold Z}).$$
\endproclaim

\demo{Proof}  The idea of the proof for the existence of the pair  
$(m, MW(f)) = (69, {\bold Z}/4{\bold Z})$  is as follows.
By \cite{MP3}, $s = (0,1,1,0,1,3)$  is the generator of  
$MW(f) = {\bold Z}/4{\bold Z}$.  As in the proof of Lemma \copy\lemmaii,
the minimal resolution  $Y$  of  $X/\langle 2s \rangle$  is
of type  $m = 104$.  The detailed proof of the existence is 
given below.

\smallbreak
Let  $g : Y \rightarrow {\bold P}^1$  be of type  $m = 104$.
By the Table in \cite{MP3}, 
$MW(g) = {\bold Z}/2{\bold Z} \times {\bold Z}/2{\bold Z}$
and we may assume that  $g$  has singular fibres
$\sum_{j=0}^{n_i - 1} \theta(i)_j$ ($i = 1, \dots, 6$)
of type  $I_{n_i}$, and two 2-torsion sections
$t_1 = (0,2,2,0,1,3), t_2 = (1,2,2,3,0,0)$, 
after suitably indexing singular fibers
so that  $[n_1, \dots, n_6] = [2,4,4,6,2,6]$.
It is easy to check the following relation
(cf. \cite{S} Lemma 8.1 or \cite{M} Formula (2..5)), where  $O_1, F$  are
the zero section and a general fiber of  $g$,
$$2t_2 \sim 2(O_1 + 2F) - (\theta(1)_1 + 
\theta(2)_1 + 2\theta(2)_2 + \theta(2)_3 +
\theta(3)_1 + 2\theta(3)_2 + \theta(3)_3 +$$
$$\theta(4)_1 + 2\theta(4)_2 + 3\theta(4)_3 +
2\theta(4)_4 + \theta(4)_5.$$
Hence we get a relation
$$D = \theta(1)_1 + \theta(2)_1 + \theta(2)_3 +
\theta(3)_1 + \theta(3)_3 +
\theta(4)_1 + \theta(4)_3 +
+ \theta(4)_5 \sim 2L$$
for some integral divisor  $L$.
Let  $\pi : {\widetilde X} \rightarrow Y$  be the
${\bold Z}/2{\bold Z}$-cover, branched along  $D$  and induced
from the above relation.  Then  $g$  induces
an elliptic fibration  $f : {\widetilde X} \rightarrow {\bold P}^1$  
so that the relatively minimal model  $(X, f)$   of
$({\widetilde X}, f)$  is of type  $m = 69$.
The inverse on  $X$  of  $O_1$  is a disjoint union of
two sections, one of which will be fixed as  $O$  of  $f$.
Now the inverse on  $X$  of the 2-torsion section  $t_1$
on  $Y$  is a disjoint union of two 4-torsion sections of
$f$.  Hence  $MW(f) = {\bold Z}/4{\bold Z}$
by the Table in \cite{MP3}.  This proves the existence
of the pair  $(m, MW(f)) = (69, {\bold Z}/4{\bold Z})$.

\smallbreak
The existence of other pairs is similar.
Here we just show which  $Y$  and  $t_1, t_2$  we should choose.
To be precise, we let
$g : Y \rightarrow {\bold P}^1$  be of type $m = 52$
(resp. $m = 97$; $m = 91$; $m = 110$; $m = 97$; $m = 104$)
and hence have singular fibers
of type  $I_{n_1} + \dots + I_{n_6}$  with  $[n_1, \dots, n_6] =
[2,1,1,6,8,6]$  (resp. $[2,6,8,2,2,4]$;
$[3,3,6,6,1,5]$; $[3,3,6,6,3,3]$; $[2,2,4,8,2,6]$; 
$[2,2,6,6,4,4]$)
and we let  $t_1 = O_1$  be the zero section 
and  $t_2 = (1,0,0,3,4,0)$  the 2-torsion section
(resp. $t_1 = (0,0,4,1,1,2)$  and  $t_2 = (1,3,4,0,0,0)$  two 2-torsion
sections;
$t_1 = O_1$  and  $t_2 = (1,1,2,2,0,0)$  a 3-torsion section;
$t_1 = O_1$  and  $t_2 = (1,1,2,2,0,0)$  a 3-torsion section;
$t_1 = (0,0,0,4,1,3)$  and  $t_2 = (1,1,2,4,0,0)$  two 2-torsion sections;
$t_1 = O_1$  and  $t_2 = (1,1,3,3,0,0)$  a 2-torsion section).
Then as in the above paragraph,
the minimal model  $X$  of a  ${\bold Z}/n{\bold Z}$-cover 
with  $n = 2$  (resp.  $n = 2$; $n = 3$; $n = 3$; $n = 2$; $n = 2$)
of  $Y$  has
an elliptic fibration  $f : X \rightarrow {\bold P}^1$,
induced from  $g$,
of type  $m = 69$  (resp. $m = 92$; $m = 32$; $m = 37$; $m = 44$; $m = 55$)
such that the inverse on  $X$
of  $t_1$  is a disjoint union of  $O$  and  $s = (0,1,1,0,0,6)$
(resp. a disjoint union of two 4-torsion sections;
a disjoint union of  $O$  and two 3-torsion sections;
a disjoint union of  $O$  and two 3-torsion sections;
a disjoint union of two 4-torsion sections;
a disjoint union of  $O$  and a  2-torsion section),
whence  $MW(f)$  is equal to   ${\bold Z}/2{\bold Z} = \{O, s\}$
(resp. ${\bold Z}/4{\bold Z}$;
${\bold Z}/3{\bold Z}$;
${\bold Z}/3{\bold Z}$; 
${\bold Z}/4{\bold Z}$;
${\bold Z}/2{\bold Z}$)
by the Table in \cite{MP3}.  

\smallbreak 
This completes the proof of the lemma and also that of
Theorem \copy\thmi.
\enddemo

\finparrafo

\head\sec Uniqueness for some of extremal elliptic $K3$ surfaces
\endhead

The goal of this section is to prove Theorem \copy\thmii. 

In the case where $MW(f) \supseteq \ZZ/2\ZZ \times \ZZ/2\ZZ$, 
namely, $m = 94, 97, 98, 103, 104, 112$, the uniqueness problem has
already been considered in \S 7 \cite{MP3} by using double sextics, 
and they are all unique.
Hence we need to  prove the cases $MW(f) \cong \ZZ/4\ZZ$, $\ZZ/5\ZZ$,
 $\ZZ/6\ZZ, \ZZ/7\ZZ$,
$\ZZ/3\ZZ\times\ZZ/3\ZZ$.

As we have seen in \S1, if $MW(f)$ has a element of order $N \ge 3$, then
$f : X \to \PP^1$ is obtained as the pull-back surface of the rational
elliptic surface, $\psi_{1,N} : E_1(N) \to X_1(N)$,
 by some morphism $g : \PP^1 \to X_1(N)$. Note that
$X_1(N)$ should be isomorphic to $\PP^1$ in our case, and this gives
a restriction on $N$.

Our proof of Theorem \copy\thmii\ consists of several steps depending on $N$.

\subhead 1.-The Case $MW(f) \cong \ZZ/4\ZZ$
\endsubhead
\bigbreak

There are 5 cases: $m = 4, 31, 44, 69, 92$.

The degree of the $j$-invariant of $E_1(4)$ is $6$,  as it has
three singular fibers $I_1^*$, $I_4$ and $I_1$. With a suitable affine
coordinate of $X_1(4)$, we may assume that these singular fibers are
over $0$, $1$ and $\infty$ respectively. Since the degree of the $j$-invariant
of $f : X \to \PP^1$ is 24, the degree of $g$ is $4$, and $g$ is branched only
at $0$, $1$ and $\infty$. By Table 7.1 in \cite{MP1} and the Riemann-Hurwitz
formula
for $g : \PP^1 \to X_1(4)$, we have the following table on
the ramification types over each branch point.
\bigbreak

\newbox\tablei
\global\setbox\tablei=\etiqueta
\centerline{\bf Table\cita}

$$
\vbox{\offinterlineskip
\halign{\tv\ #&&\tv\centro{#}&\tv#\cr
\hline
$m$ &  The ramification types over $0$, $1$ and $\infty$  &\cr \hline
$4$ & $(4)$, $(4)$, $(1, 1, 1, 1)$ &\cr \hline
$31$ & $(2, 2)$, $(4)$, $(2, 1, 1)$ &\cr \hline
$44$ & $(4)$, $(3, 1)$, $(2, 1, 1)$ &\cr \hline
$69$ & $(2, 2)$, $(3, 1)$, $(3, 1)$  &\cr \hline
$92$ & $(4)$, $(2, 1, 1)$, $(3, 1)$ &\cr \hline}
}
$$

Here the notation $(e_1,...,e_k)$ means that $g^{-1}(p)$ 
($p \in \{ 0, 1, \infty \}$) consists of $k$ points, $q_1$,...,$q_k$, and
the ramification index at $q_j$ is $e_j$.

To show the uniqueness of surfaces, it is enough to
show that $g$ assigned with the ramification types as above is unique
up to covering isomorphisms over $X_1(4)$. For this purpose, the following lemma
is
important.

\newbox\gpx
\global\setbox\gpx=\etiqueta
\proclaim{Lemma\cita} Let $g : \PP^1 \to X_1(4)$ be a degree 4 map in 
Table \copy\tablei. Let $\alpha : C \to \PP^1$ be the Galois closure, and put
$\hat g = g \circ \alpha$. Then we have the following:

\smallbreak\item{}$m = 4$: $g = \hat g$ and $g$ is a $4$-fold cyclic covering.

\smallbreak\item{}$m = 31$: $\deg \hat g = 8$, $C \cong \PP^1$ and
$\Gal(\hat{g}) \cong {\Cal D}_8$.

\smallbreak\item{}$m = 44, 92$: $\deg \hat g = 24$, $C \cong \PP^1$ and $\Gal
(\hat{g}) \cong
\SS_4$. 

\smallbreak\item{}$m = 69$: $\deg \hat g = 12$, $C \cong \PP^1$ and $\Gal
(\hat{g}) \cong {\Cal A}_4$.
\endproclaim

\demo{Proof} The monodromy around the branch points gives a permutation 
representation of $\pi_1(\PP^1 \setminus \{0, 1, \infty\})$ to $\SS_4$; the
basic loops $\gamma_0$, $\gamma_1$ and $\gamma_{\infty}$ about $0$, $1$ and
$\infty$, respectively map to permutations $\sigma_0$, $\sigma_1$ and 
$\sigma_{\infty}$. The cycle structure of each permutation is the
same as the ramification type over the corresponding point.
These permutations satisfy the identity 
$\sigma_0\sigma_1\sigma_{\infty} = 1$ in $\SS_4$ and generate a transitive
subgroup, $G$, in $\SS_4$. 
Note that this $G$ is nothing but the Galois group of
$\hat g : C \to X_1(4)$. We apply this argument to each case, and obtain
the following table:
\bigbreak
\vfill\eject
\newbox\tableii
\global\setbox\tableii=\etiqueta

\centerline{\bf Table\cita}

$$
\vbox{\offinterlineskip
\halign{\tv\ #&&\tv\centro{#}&\tv\ #& \tv#\cr
\hline
$m$ &  The cycle structure of $\sigma_0$, $\sigma_1$ and $\sigma_{\infty}$
& $G$   &\cr \hline
$4$ & $(4)$, $(4)$, $(1, 1, 1, 1)$ & $\ZZ/4\ZZ$ &\cr \hline
$31$ & $(2, 2)$, $(4)$, $(2, 1, 1)$ & ${\Cal D}_8$ &\cr \hline
$44$ & $(4)$, $(3, 1)$, $(2, 1, 1)$ & $\SS_4$ &\cr \hline
$69$ & $(2, 2)$, $(3, 1)$, $(3, 1)$ & ${\Cal A}_4$ &\cr \hline
$92$ & $(4)$, $(2, 1, 1)$, $(3, 1)$ & $\SS_4$ &\cr \hline
}}
$$

Now all we need to show are the assertions: $C \cong \PP^1$. 
Our argument is based
on the following elementary fact:

\newbox\facti
\global\setbox\facti=\etiqueta
\proclaim{Fact\cita} Let $x$ be a point on $C$, and put
 $G_x = \{ \tau \in G | \tau(x) = x \}$. Then 
$$
\vbox{\offinterlineskip
\halign{\tv\ #&&\tv\centro{#}&\tv#\cr
\hline
$G$ & The order of $G_x$ &\cr \hline
$\ZZ/4\ZZ$ & $1, 2, 3$ &\cr \hline
$\SS_4$ & $1, 2, 3, 4$  &\cr \hline
$\AA_4$ & $1, 2, 3$  &\cr \hline
${\Cal D}_8$ & $1, 2, 4$ &\cr \hline
}}
$$
\endproclaim

We prove $C \cong \PP^1$ case by case.

$m = 4$: As $G = \ZZ/4\ZZ$, $\deg \hat g = \deg g$, and $\alpha$ is the
identity.

$m = 31$: Since $G = {\Cal D}_8$, $\deg \alpha = 2$. Let $\iota$ be an element
 of order $2$ such that $C/\langle \iota \rangle \cong \PP^1$. As $g$ is 
not Galois, $\iota \not\in \text{center of }{\Cal D}_8$. If $\alpha$ is
branched over $g^{-1}(0)$, then ${\hat g}^{-1}(0)$ consists of two points,
each of which has the ramification index $4$. This means that
$\iota$ belongs to the center of ${\Cal D}_8$, which leads us to a contradiction.
Hence the branch points of $\alpha$ are two points in $g^{-1}(\infty)$ which
are unramified points of $g$. Hence $C \cong \PP^1$.

$m = 44,\, 92$: By Fact \copy\facti\ and $\Gal (C/\PP^1) \cong \SS_4$, points
over 
$0$, $1$ and $\infty$ have the ramification indices $4$, $3$ and $2$, 
respectively. By the Riemann-Hurwitz formula, we have $C \cong \PP^1$.

$m = 69$: By Fact \copy\facti, points over $0$, $1$ and $\infty$ have 
the ramification indices $2$, $3$ and $3$, respectively. By
the Riemann-Hurwitz formula, $C \cong \PP^1$.

This completes our proof for Lemma \copy\gpx.
\enddemo

The following classical fact is a key to prove Theorem \copy\thmii\
in the case where
$MW(f) \cong \ZZ/4\ZZ$.

\newbox\factii
\global\setbox\factii=\etiqueta
\proclaim{Fact\cita ([Na] pp.31 -32)} 
For a suitable choice of an affine coordinate, $w$ and $z$, 
of $X_1(4)$ and $\PP^1$, respectively,
the map in Table \copy\tableii\ can be given by the rational 
functions as follows:

$$
\matrix
w = z^4 & m = 4 \\
w = - \frac {(z^4 -1)^2}{4z^2} & m = 31 \\
w = \left (\frac {z^4 + 2\sqrt 3 z^2 - 1}{z^4 - 2\sqrt 3 z^2 -1} \right )^3
& m = 69 \\
w = \frac {(z^8 + 14 z^4 + 1)^3}{108z^4(z^4 - 1)^4} & m = 44, 92
\endmatrix
$$
\endproclaim
\

Fact \copy\factii\ 
implies that the Galois coverings described in Lemma \copy\gpx\ are 
essentially unique up to isomorphisms over $\PP^1$. 
The morphisms $g$ in Lemma \copy\gpx\
are corresponding to subgroups of index $4$ of $G$, and for each case, such
subgroups are conjugate to each other. This shows that the pull-back
morphisms, $g$, are unique up to covering isomorphisms over $X_1(4)$. Therefore
we have Theorem \copy\thmii\ in the case where $MW(f) \cong \ZZ/4\ZZ$.

\remark{Remark\cita} We can prove the uniqueness for $m = 94, 98, 103, 112$ in
a similar way to the case $MW(f) \cong \ZZ/4\ZZ$.
\endremark

\subhead 2.-The Case $MW(f) \cong \ZZ/5\ZZ$
\endsubhead
\bigbreak

There are 3 cases: $m = 10, 49, 105$.
 $f : X \to \PP^1$ is obtained as the pull-back surface of 
 $\psi_{1,5} : E_1(5) \to X_1(5)$ by a degree $2$ map 
 $g : \PP^1 \to X_1(5)$. There are four singular fibers for
 $\psi_{1,5}$, which are $I_5$, $I_5$, $I_1$, $I_1$. By \cite{MP1} Table 5.3,
$E_1(5)$ is given by the following Weierstrass equation:
$$
y^2  =  x^3  - 3(s^4 - 12 s^3 + 14s^2 + 12 s + 1)x 
  + 2(s^6 - 18 s^5 + 75 s^4 + 75 s^2 + 18s + 1),
$$
where $s$ is an affine coordinate of $X_1(5) \cong \PP^1$. 
The two $I_5$ fibers are over 
$s = 1$ and
$s = \infty$, and the two $I_1$ fibers are over $s = (11 \pm 5\sqrt 5)/2$.

$m = 10$: The pull-back morphism $g$ is branched at 
$s = 0$ and $\infty$, and such a morphism is unique.

$m = 49$: There are $4$ possible cases for the pull-back morphism 
depending on the branch points as follows:
$$
\vbox{\offinterlineskip
\halign{\tv\ #&&\tv\centro{#}&\tv#\cr
\hline
  & The branch points of $g$ &\cr \hline
(1)  & $0$ and $(11 + 5\sqrt 5)/2$  &\cr \hline
(2)  & $0$ and $(11 - 5\sqrt 5)/2$  &\cr \hline
(3)  & $\infty$ and $(11 + 5\sqrt 5)/2$  &\cr \hline
(4)  & $\infty$ and $(11 - 5\sqrt 5)/2$  &\cr \hline
}}
$$

\proclaim{Proposition\cita} 
There exists $\varphi$ in Question \copy\cuestion\ between
the two pull-back surfaces for either $(1)$ and $(4)$, or $(2)$ and $(3)$, 
while there is no such $\varphi$ between the two pull-back surfaces for
other combinations.
\endproclaim

\demo{Proof}  Consider an automorphism, $\tau$, of $E_1(5) \to
X_1(5)$ given by
$$
\tau :  (x, y, s) \mapsto \left ( \frac 1{s^2}x, \frac 1{s^3}y, -\frac 1s
\right).
$$
With $\tau$, the points $0$ and $(11 + 5\sqrt 5)/2$ map to
$\infty$ and $(11 - 5\sqrt 5)/2$, respectively. Our first assertion
follows from this fact. For the second, by using $\tau$,
it is enough to show that
there is no $\varphi$ in Question 
\copy\cuestion\ between the pull-back surfaces for
$(1)$ and $(2)$.

Let $f_i : X_i \to \PP^1$ $(i = 1, 2)$ be the pull-back surfaces for $(1)$
 and $(2)$, respectively. Suppose that there exists $\varphi : X_1 \to X_2$
as Question \copy\cuestion. Then we have

\proclaim{Claim\cita} $\varphi$ induces an automorphism $\tilde {\varphi} :
X_1(5) \to X_1(5)$ such that $0 \mapsto \infty$,
$\infty \mapsto 0$, 
$(11 + 5\sqrt 5)/2 \mapsto (11 - 5\sqrt 5)/2$, and 
$(11 - 5\sqrt 5)/2 \mapsto (11 + 5\sqrt 5)/2$.
\endproclaim

Since there is no fractional linear transformation as above, the second
assertion follows.

\demo{Proof of the Claim}  Let $\iota_i$ $(i = 1, 2)$ be fiber
preserving involutions
on $X_i$ $(i = 1, 2)$ determined by the pull-back morphisms $g_i$. Let 
$\bar {\varphi}$ and $\bar {\iota}_i$ $(i = 1, 2)$ be the restrictions of
each morphism to the
zero sections of $X_1$ and $X_2$. 
$\varphi^{-1}\circ\iota_2\circ\varphi$ gives rise to 
another fiber preserving involution on $X_1$. 
With $\varphi^{-1}\circ\iota_2\circ\varphi$, $I_{10}$, $I_5$, $I_2$ fibers
map to $I_{10}$, $I_5$, $I_2$ fibers, respectively. Hence
$\bar {\varphi}^{-1}\circ
\bar {\iota}_2\circ \bar {\varphi} = \bar {\iota}_1$ or $id$, but the latter
case does not occur since $\iota_2 \neq id$. Thus we have an isomorphism
$\tilde {\varphi} : X_1(5) \to X_1(5)$, and it is easy to
see that $\tilde {\varphi}$ has the desired property.
\enddemo

$m = 105$: Likewise $m = 10$, the surface is unique.
\enddemo

\subhead 3.-The case $MW(f) \cong \ZZ/6\ZZ$
\endsubhead
\bigbreak

There are $5$ cases: $m = 35, 53, 63, 95, 108$. 
For all cases,
 $f : X \to \PP^1$ is obtained as the pull-back surface of 
 $\psi_{1, 6} : E_1(6) \to X_1(6)$ by a degree $2$ map 
 $g : \PP^1 \to X_1(6)$, and they are unique by a similar argument
 to $m = 10$.

\subhead 4.-The case $MW(f) \cong \ZZ/7\ZZ$
\endsubhead
\bigbreak

There is only one case: $m = 30$.
 In this case, $f : X \to \PP^1$ is obtained as the pull-back surface of
 $\psi_{1,7} : E_1(7) \to \PP^1$. Comparing the $j$-functions of both surfaces, 
 we know that the degree of the pull-back morphism is $1$, i.e.,
 $X$ is isomorphic to $E_1(7)$. This implies the uniqueness.

\subhead 5.-The case $MW(f) \cong \ZZ/3\ZZ\times \ZZ/3\ZZ$
\endsubhead
\bigbreak

There is only one case: $m = 110$.
$f : X \to \PP^1$ is obtained as the pull-back surface of 
$\psi_{3, 3} : E_3(3) \to X_3(3)$ by a degree $2$ map 
$g : \PP^1 \to X_3(3)$. There are four singular fibers for
$\psi_{3,3}$, which are all $I_3$. By \cite{MP1} Table 5.3,
$E_3(3)$ is given by the following Weierstrass equation:
$$
y^2  =  x^3  + (- 3s^4 + 24s)x + (2s^6 + 40 s^3 - 16)
$$
where $s$ is an inhomogeneous coordinate of $X_3(3) \cong \PP^1$. 
The four $I_3$ fibers are over
$-1$, $-\omega$, $-\omega^2$ and $\infty$, where $\omega = \exp(2\pi\sqrt
{-1}/3)$.

Consider two fiber preserving automorphisms of $E_3(3)$:
$$
\tau_1: (x, y, s) \mapsto 
\left ( \frac {-3}{(s + 1)^2}x, \frac {3\sqrt {-3}}{(s + 1)^3}y,  
\frac {- s + 2}{s + 1} \right),
$$
and 
$$
\tau_2: (x, y, s) \mapsto  (\omega x, y, \omega s).
$$
These automorphisms induce permutations of the $I_3$ fibers. Since $X$ is a
double
covering of $E_3(3)$,
 it is uniquely determined by the branch locus which is
two $I_3$ fibers.
Therefore,
using $\tau_1$ and $\tau_2$, we can show that $f : X \to \PP^1$ is unique.

Summing up, we have Theorem \copy\thmii.  

\finparrafo

\enddocument